# OPERATOR NORM CONSISTENT ESTIMATION OF LARGE-DIMENSIONAL SPARSE COVARIANCE MATRICES[1]

### By Noureddine El Karoui

#### *University of California, Berkeley*

Estimating covariance matrices is a problem of fundamental importance in multivariate statistics. In practice it is increasingly frequent to work with data matrices $X$ of dimension $n \times p$, where $p$ and $n$ are both large. Results from random matrix theory show very clearly that in this setting, standard estimators like the sample covariance matrix perform in general very poorly.

In this "large $n$, large $p$" setting, it is sometimes the case that practitioners are willing to assume that many elements of the population covariance matrix are equal to 0, and hence this matrix is sparse. We develop an estimator to handle this situation. The estimator is shown to be consistent in operator norm, when, for instance, we have $p \asymp n$ as $n \to \infty$. In other words the largest singular value of the difference between the estimator and the population covariance matrix goes to zero. This implies consistency of all the eigenvalues and consistency of eigenspaces associated to isolated eigenvalues.

We also propose a notion of sparsity for matrices, that is, "compatible" with spectral analysis and is independent of the ordering of the variables.

**1. Introduction.** Estimating covariance matrices is the cornerstone of much of multivariate statistics. Theoretical contributions (see [2], Chapter 7, [14, 18]) have been highlighting for a long time the fact that for various loss functions, one could improve on the sample covariance matrix as an estimator of the population covariance matrix, in a nonasymptotic setting.

The "large $n$, large $p$" context, that is, multivariate analysis of datasets for which both the number of observations, $n$, and the number of variables,

Received July 2007.

[1]Supported in part by NSF Grant DMS-06-05169 and by SAMSI.

$p$, are large, is, in a somewhat different setting, highlighting the deficiency of this estimator. To be more precise, when we refer to "large $n$, large $p$" problems, we generally mean that $p \asymp n$, that is, $p = \mathrm{O}(n)$ and $n = \mathrm{O}(p)$, and $p \to \infty$. If $p/n$ has a nonzero limit as $n \to \infty$, results from random matrix theory [21] make clear that in this asymptotic setting, even at just the level of eigenvalues, the sample covariance matrix will not lead to a consistent estimator. We refer to [12] for a more thorough introduction to these ideas and the consequences of the results for statistical practice.

This is naturally very problematic since this class of results suggests that the sample covariance matrix contains little reliable information about the population covariance. This realization has helped generate a significant amount of work in mathematics, probability and theoretical statistics and the behavior of many hard to analyze quantities is now quite well understood. For instance, under strong distributional assumptions, one can characterize the fluctuation behavior of the largest eigenvalue of sample covariance matrices for quite a large class of population covariance (see, e.g., [11] for recent results), or the fluctuation behavior of linear functionals of eigenvalues (see [1, 3, 19]). However, until very recently there has been less work in the direction of using these powerful results for the sake of concrete data analysis.

Of course, since this inconsistency phenomenon is now fairly well-known, remedies have been proposed. For instance, the interesting paper [20] proposes to shrink the sample covariance matrix toward the identity matrix using a shrinkage parameter chosen from the data. In [12], a nonparametric estimator of the spectrum is proposed and shown to be consistent in the sense of weak convergence of distributions. The method in [12] uses convex optimization, random matrix theory (the generalization of [21] found in [22]) and ideas from nonparametric function estimation. These estimation methods rely on asymptotic properties of eigenvalues, and as a starting point for estimation of the full covariance matrix, they are essentially trying to get an estimator, that is, equivariant under the action of the orthogonal (or unitary) group. In other words, the "basis" in which the data are given is not taken advantage of, and the premise of such an analysis is that we should be able to find good estimators in any "basis." While ideally researchers will be able to come up with strategies to solve the estimation problem at this level of generality, it is reasonable to expect that taking advantage of the representation of the data we are given should or might help finding good estimators.

In particular, it is often the case that data analysts are willing to assume that the basis in which the data are given is somewhat nice. Often this translates into assumption that the population covariance matrix has a particular structure in this basis, which should naturally be taken advantage of. In this situation, it becomes natural to perform certain forms of regularization by working directly on the entries of the sample covariance matrix. Various strategies have been proposed (see [4, 17]) that try to take advantage of the assumed structure. The very interesting paper [7] proposed the idea of "banding" covariance matrices when it is known that the population covariance has small entries far away from the diagonal. The idea is to put to zero all coefficients that are too far away from the diagonal and to keep the other ones unchanged. Remarkably, in [7], the authors show consistency of their estimator in spectral (a.k.a. operator) norm, a very nice result. In other words, they show that the largest singular value of the difference between their estimator and the population covariance matrix goes to zero as both dimensions of the matrices go to infinity and, for instance, when $p/n$ has a finite limit. The requirement of estimating consistently in spectral norm is a very interesting idea (which we adopt in this paper), since then one can deduce easily many results concerning consistency of eigenvalues and eigenspaces. We make this remark more precise in Section 3.5, using different arguments than those used by Bickel and Levina in [7].

It might be argued that ideas such as banding essentially assume that one knows a "good" ordering of the variables. As a matter of fact, if we start with a matrix with entries small or zero away from the diagonal and reorder the variables, the new covariance matrix we obtain may not have only small entries away from the diagonal. In some situations, for instance, time series analysis, the order of the variables has a statistical/scientific meaning and so using it makes sense. However, in many data-analytic problems, there is no canonical ordering of the variables.

Hence to tackle these situations, a natural requirement is to find an estimator which is equivariant under permutations of the variables. We call such estimators permutation-equivariant. Such an estimator would take advantage of the particular nature of the basis in which the data are given, while not requiring the user to find a permutation of the order of the variables that makes the analysis particularly simple. Searching for such a permutation would—in general—be practically infeasible. Note that regularizing the estimator by applying the same function to each of the entries of the matrix leads to permutation-equivariant estimators.

A subject of particular practical interest is the estimation of sparse covariance matrices (see, e.g., [9]) because they are appealing to practitioners for several reasons, including interpretability, presumably ease of estimation and the practically often encountered situation where while many variables are present in the problem, most of them are correlated to only "a few" others.

In this paper we propose to estimate sparse matrices by hard thresholding small entries of the sample covariance matrix and putting them to zero. We propose a combinatorial view of the problem, inspired in part by classical ideas in random matrix theory, going back to [25]. The notion of sparsity we propose is flexible enough that it makes the proofs manageable and at the same time rich enough that it captures many practically natural situations.

We show that our estimators are consistent in spectral norm, in the case of both the sample covariance and the sample correlation matrix. No assumptions of normality of the data are required, only the existence of certain moments. As is to be expected, the larger the number of moments available, the easier the task and the larger the class of matrices we can estimate consistently.

Finally, we note that at the same time as we were researching these questions and independently, similar questions were approached from a very different point of view by [8].

## 2. Sparse matrices: concepts and definitions.

One conceptual difficulty of this problem is to define a notion of sparsity for matrices that is compatible with spectral analysis. Just as in the case of norms, extending straightforwardly the notions from vectors to matrices can be somewhat unhelpful. In the norm case, the Frobenius norm—the extension of the $\ell_2$ (vector) norm to matrices—is, for instance, known to not be as good as other matrix norms from a spectral point of view. Similarly here, we will explain that just counting the number of 0's in the matrix—the canonical sparsity notion for vectors—does not yield a "good" notion of sparsity when one investigates the spectral properties of matrices.

Let us illustrate our problem on a concrete example. Consider now two $p \times p$ symmetric matrices with the same number of nonzero coefficients:

$$E_1 = \begin{pmatrix} 1 & \frac{1}{\sqrt{p}} & \frac{1}{\sqrt{p}} & \cdots & \frac{1}{\sqrt{p}} \\ \frac{1}{\sqrt{p}} & 1 & 0 & \cdots & 0 \\ \vdots & \vdots & \ddots & \ddots & 0 \\ \frac{1}{\sqrt{p}} & 0 & 0 & 1 & 0 \\ \frac{1}{\sqrt{p}} & 0 & 0 & \cdots & 1 \end{pmatrix}$$

and

$$E_2 = \begin{pmatrix} 1 & \dfrac{1}{\sqrt{p}} & 0 & \dots & 0 \\ \dfrac{1}{\sqrt{p}} & 1 & \dfrac{1}{\sqrt{p}} & \dots & 0 \\ 0 & \dfrac{1}{\sqrt{p}} & 1 & \dfrac{1}{\sqrt{p}} & \dots \\ \vdots & \ddots & \ddots & \ddots & \dfrac{1}{\sqrt{p}} \\ 0 & \dots & 0 & \dfrac{1}{\sqrt{p}} & 1 \end{pmatrix}.$$

Using the Schur decomposition of $E_1$ to compute its characteristic polynomial (see also Section 3.3), we see easily that its eigenvalues are $(p-2)$ 1's and $1 + \sqrt{p-1}/\sqrt{p}$ and $1 - \sqrt{p-1}/\sqrt{p}$. On the other hand, $E_2$ is a well-known matrix, for instance, in numerical analysis, and its eigenvalues are $\{1 + 2\cos(k\pi/(p+1))/\sqrt{p}\}_{k=1}^{p}$. Hence, the extreme eigenvalues of these matrices are very different, but they have the same number of nonzero coefficients and their nonzero coefficients have the same values. This raises the question of trying to come up with an alternative notion of sparsity that, while encompassing the canonical notion of "having a large number of zeros," might be better suited for the study and the understanding of spectral properties of matrices.

2.1. *Matrix sparsity: proposed definition.* To describe our proposal, we need to introduce several concepts from graph theory and combinatorics. For the sake of readability we detail them here; they can also be found in, for instance, [23], Section 4.7. To each population covariance matrix, $\Sigma_p$, it is natural to associate an adjacency matrix $A_p(\Sigma_p)$, in the following fashion:

$$A_p(i,j) = 1_{\sigma(i,j)\neq 0}.$$

This matrix $A_p$ can in turn be viewed as the adjacency matrix of a graph $\mathcal{G}_p$, with $p$ vertices, corresponding to the variables in our statistical problem. We call a walk on this graph closed if it starts and finishes at the same vertex. The length of a walk is the number of edges it traverses. By definition, we call

$$\mathcal{C}_p(k) = \{\text{closed walks of length } k \text{ on the graph with adjacency matrix } A_p\}$$

and

$$\phi_p(k) = \mathrm{Card}\{\mathcal{C}_p(k)\}.$$

Note that we have

$$\phi_p(k) = \mathrm{trace}(A_p^k).$$

The following two figures show the graphs corresponding to the adjacency matrices of $E_1$ and $E_2$:

Graph corresponding to $E_1$:                    Graph corresponding to $E_2$:

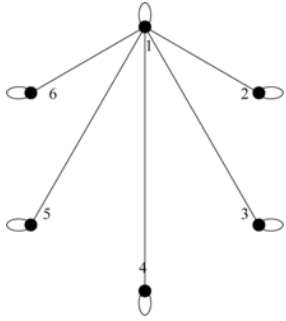          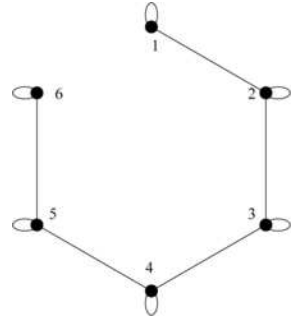

DEFINITION 1.   We say that a sequence of covariance matrices $\{\Sigma_p\}_{p=1}^{\infty}$ is $\beta$-sparse if the graphs associated to them via $A_p$'s have the property that

$$\forall k \in 2\mathbb{N} \qquad \phi_p(k) \leq f(k) p^{\beta(k-1)+1},$$

where $f(k) \in \mathbb{R}^+$ is independent of $p$ and $0 \leq \beta \leq 1$.

We say that a sequence of matrices is asymptotically $\beta$-sparse if it is $\beta + \varepsilon$ sparse for any $\varepsilon > 0$.

We call $\beta$ an index of sparsity of the sequence of matrices.

For short, we say that a matrix is $\beta$-sparse instead of saying that a sequence of matrices is $\beta$-sparse when this shortcut does not cause any confusion.

Here are a few simple examples of matrices that are sparse according to our definition:

1. **Diagonal matrices.** In the case of diagonal matrices, $A_p = \mathrm{Id}_p$, and $\mathcal{G}_p$ consists only of self-loops at each vertex. Hence $\phi(k) = p$, for all $k$. So a diagonal matrix is 0-sparse.

2. **Matrices with at most $M$ nonzero elements on each line.** For these matrices, the corresponding $\mathcal{G}_p$ has at most $M$ edges at each vertex. A simple enumeration shows that $\phi(k) \leq pM^{k-1}$. So these matrices are also 0-sparse.

3. **Matrices with at most $Mp^{\alpha}$ nonzero elements on each line.** The same argument shows that $\phi(k) \leq p(Mp)^{\alpha(k-1)}$. So these matrices are $\alpha$-sparse. In particular, full matrices are 1-sparse.

4. **Matrices with at most $M(\log p)^r$ nonzero elements on each line.** We have, by simple counting arguments, $\phi_p(k) \leq pM^{k-1}(\log p)^{r(k-1)}$. These matrices are therefore $\beta$-sparse for any $\beta > 0$ and asymptotically 0-sparse.

Given a matrix $S_p$, we can associate to the corresponding $\mathcal{G}_p$ a set of weights on the edges, by simply setting the weight of the edge joining vertices $i$ and $j$ to $S_p(i, j)$. Similarly, for a walk, we have

DEFINITION 2 (Weight of a walk). Given $\gamma$, a closed walk of length $k$: $\gamma: i_1 \to i_2 \to i_3 \to \cdots \to i_k \to i_{k+1} = i_1$, and a matrix $S_p$, we call $w_\gamma$ the weight of the walk $\gamma$. By definition it is

$$w_\gamma = S_p(i_1, i_2) S_p(i_2, i_3) \cdots S_p(i_k, i_1).$$

We conclude this section by the following simple but important remark:

$$\text{trace}(S_p^k) = \sum_{\gamma \in \mathcal{C}_p(k)} w_\gamma.$$

2.2. *Remarks on the notion of sparsity proposed.* It is clear that if we change the order of the variables in our statistical problem, we do not change the "index of sparsity" of our matrices. This is essentially obvious from the graph representation of the problem. From a more algebraic standpoint, if the permutation that is applied is encoded as a permutation matrix $P$, the covariance in the permuted problem is simply $P'\Sigma_p P$ and the new adjacency matrix is $P'A_p P$ (this matrix is indeed an adjacency matrix). Since $P'P = \text{Id}_p$, we have $\text{trace}((P'A_p P)^k) = \text{trace}(A_p^k)$, and hence the sparsity index is unchanged when we permute the variables.

We also note that we could replace the notion of $\beta$-sparsity we use by the requirement that, for some $C > 0$,

$$\phi_p(k) \leq C^k p^{1+\beta k} \qquad \forall k \in 2\mathbb{N},$$

which is equivalent, if $\|\!\|\cdot\|\!\|_2$ represents the operator norm or largest singular value of a matrix, to

$$\|\!\|A_p\|\!\|_2 \leq C p^\beta.$$

This would result in minor differences in the theorems that follow and might be slightly simpler to apply when the only information available concerns the largest eigenvalue of $A_p^2$. From a combinatorial point of view, the notion we use in this paper is a bit more natural and this is what directed our choice.

Finally, we also note that we could replace our notion of sparsity by

$$\forall k \leq k_0, \ k \in 2\mathbb{N}, \qquad \phi_p(k) \leq f(k) p^{\beta(k-1)+1},$$

and call the matrices having this property $\beta$-sparse at order $k_0$. The proofs below would still be valid provided $k_0$ is large enough, the minimum required size for $k_0$ being related to the number of moments of the random variables making up our data matrix.

Let us now return to the notion of $\beta$-sparsity proposed in Definition 1. It is clear that the smaller $\beta$ is, the sparser the matrix. In particular, if $\beta_0 \leq \beta_1$, a matrix which is $\beta_0$-sparse is also $\beta_1$-sparse. As we will shortly show, the class of $\beta$-sparse matrices is stable by addition, which implies that the sum of a $\beta_0$-sparse and a $\beta_1$-sparse matrix is $(\beta_0 \vee \beta_1)$-sparse.

We conclude this discussion with two properties of $\beta$-sparse matrices.

FACT 1. *The set of $\beta$-sparse matrices is stable by addition. In other words, the sum of two $\beta$-sparse matrices is $\beta$-sparse.*

PROOF. We call $B_0$ and $B_1$ our "initial" $\beta$-sparse matrices, and $B_2$ their sum. $A_2$, the adjacency matrix of $B_2$, is not the sum of $A_0 + A_1$. In particular, edges that are present in both $A_0$ and $A_1$ may not be present in $A_2$. However, if we add edges to $A_2$, we increase $\phi_p^{(2)}(k)$, the number of closed walks of length $k$ on $A_2$. So in checking the sparsity index of $B_2$, we can work with $\widetilde{A}_2$, which contains all edges in $A_0$ and $A_1$, and contains the graph corresponding to $A_2$ as a subgraph of its own graphical representation. More algebraically, the definition of $\widetilde{A}_2$ is

$$\widetilde{A}_2(i,j) = \min(A_0(i,j) + A_1(i,j), 1) = 1_{A_1(i,j)=1} + 1_{A_0(i,j)=1} 1_{A_1(i,j)=0}.$$

We can write $\widetilde{A}_2 = \widetilde{A}_0 + A_1$, with $\widetilde{A}_0(i,j) = 1_{A_0(i,j)=1} 1_{A_1(i,j)=0}$. Note that $\widetilde{A}_0$ is a symmetric adjacency matrix, may have zeroes where $A_0$ has ones, but does not have ones where $A_0$ has zeroes. So the graph corresponding to $\widetilde{A}_0$ is a subgraph of the graph corresponding to $A_0$. In particular, trace($\widetilde{A}_0^{2k}$) $\leq$ trace($A_0^{2k}$).

The matrices $\widetilde{A}_0$, $A_1$ and $\widetilde{A}_2$ are all symmetric, so when dealing with their eigenvalues we can apply standard results for symmetric matrices. Using Lidskii's theorem (see [6], Corollary III.4.2), we know that

$$\lambda^{\downarrow}(\widetilde{A}_2) \prec \lambda^{\downarrow}(\widetilde{A}_0) + \lambda^{\downarrow}(A_1),$$

where $\lambda^{\downarrow}(A_1)$ is the vector of decreasing eigenvalues of $A_1$ and the sign $\prec$ means that the left-hand side is majorized by the right-hand side (see [6], page 28 for a definition, if needed). Now the functions $h(x) = x^{2k}$ are convex and we therefore have, using standard results in the theory of majorization ([6], Theorem II.3.1),

$$\text{trace}(\widetilde{A}_2^{2k}) \leq \sum [\lambda_j(\widetilde{A}_0) + \lambda_j(A_1)]^{2k} \leq 2^{2k-1} \sum \lambda_j(\widetilde{A}_0)^{2k} + \lambda_j(A_1)^{2k}$$

$$\leq 2^{2k-1} \text{trace}(A_0^{2k} + A_1^{2k}).$$

Because $A_0$ and $A_1$ are $\beta$-sparse, we see that $\widetilde{A}_2^{2k}$ is. And because we have seen that

$$\operatorname{trace}(A_2^{2k}) \leq \operatorname{trace}(\widetilde{A}_2^{2k}),$$

we conclude that $B_2$ is $\beta$-sparse. $\square$

Fact 2. *The set of $\beta$-sparse matrices is not stable by inversion or multiplication.*

Proof. To prove this fact, it suffices to provide an example. Let us consider

$$\Sigma_p = \begin{pmatrix} 1 & \alpha & \alpha & \dots & \alpha \\ \alpha & 1 & 0 & \dots & 0 \\ \vdots & \vdots & \ddots & \ddots & 0 \\ \alpha & 0 & 0 & 1 & 0 \\ \alpha & 0 & 0 & \dots & 1 \end{pmatrix}.$$

As explained in Section 3.3 below, this matrix is a rank-2 perturbation of $\mathrm{Id}_p$, and its eigenvalues and eigenvectors are known. Note that $\Sigma_p$ is 1/2-sparse. Also, if $\alpha \leq 1/\sqrt{p-1}$, $\Sigma_p$ is a positive semidefinite matrix, and hence is a covariance matrix.

Using the expressions for eigenvalues and eigenvectors found below, we see that $\Sigma_p^2$ is full of nonzero entries, and hence is 1-sparse. As a matter of fact, it is easily checked that if $i > j \geq 2$, $\Sigma_p^2(i,j) = \alpha^2$. So there is no stability by multiplication, for otherwise $\Sigma_p^2$ would be 1/2-sparse.

Using the classic expression for inverses of low-rank perturbations of matrices found, for example, in [16], page 19, we see that $\Sigma_p^{-1}$, when it exists, is full of nonzero entries and hence is 1-sparse. As a matter of fact, it is easily checked that if $\alpha^2 \neq 1/(p-1)$, and $i > j \geq 2$, $\Sigma_p^{-1}(i,j) = -\alpha^2/(\alpha^2(p-1)-1)$. So there is no stability by inversion either. $\square$

## 3. Estimation by entrywise thresholding.
To avoid any confusion as to the meaning of the results to be proved, we remind the reader that the spectral norm of a matrix $A$ is defined (see [16], page 295) as $\|A\|_2 = \max\{\sqrt{\lambda}: \lambda$ an eigenvalue of $A^*A\}$; in other words, it is the largest singular value of $A$. When $A$ is a symmetric matrix, $\|A\|_2$ coincides with the spectral radius of $A$: $\rho(A) = \max_i |\lambda_i(A)|$. In what follows, we use interchangeably the terms spectral norm and operator norm.

When we say that we threshold a variable $x$ at level $t$ we mean that we keep (or replace $x$ by) $x 1_{|x| \geq t}$. We also refer to this operation as hard thresholding. Our final remark concerns notation: in what follows, $C$ refers to a generic constant independent of $n$ and $p$. Its value may change from display to display when there is no risk of confusion about the meaning of the statements. If there is, we also use $K$ or $C'$ and they play the same role as $C$.

3.1. *Estimation of sparse covariance or correlation matrices.* We first prove an intermediate result concerning the Gaussian MLE estimator when it is known that the mean of the data is zero (Theorem 1). This is a stepping stone to the more practically relevant results concerning the sample covariance matrix (Theorem 2) and the sample correlation matrix (Theorem 3). The proofs of these later results are essentially the same as that of Theorem 1, but the proof of Theorem 1 is technically a bit less complicated and highlights the key ideas. In Section 3.2, we explain how these results can be extended to nonsparse matrices that are approximable by sparse matrices. In Section 3.2.1, we discuss a strengthening of Theorems 1, 2 and 3, whose possibility was suggested to us by an insightful question of Professor Peter Bickel, which allows us to get rid of assumption (ii) in Theorem 1. (This strengthening is postponed to this later section for the sake of clarity.)

We refer the reader to Section 3.5 for detailed explanations of the consequences for eigenvalues and eigenvectors of Theorems 1, 2 and 3 as well as their extensions. Finally, we stress that, unless otherwise noted, all of our results are obtained when we have $p \asymp n$ as $n \to \infty$ (allowing the ratio $p/n$ to have, for instance, a finite nonzero limit), that is, in the "large $n$, large $p$" setting.

THEOREM 1. *Suppose $X$ is an $n \times p$ matrix, and assume that $p \asymp n$ as $n \to \infty$. Suppose that the rows of $X$ are independent and identically distributed and denote them by $\{X_i\}_{i=1}^n$. Call $\Sigma_p$ the covariance matrix of the vector $X_1$. We also work under the following assumptions:*

(i) *Suppose $\Sigma_p$ is $\beta$-sparse with $\beta = 1/2 - \eta$ and $\eta > 0$.*

(ii) *Suppose that the nonzero coefficients of $\Sigma_p$ are all greater in absolute value than $Cn^{-\alpha_0}$, with $0 < \alpha_0 = 1/2 - \delta_0 < 1/2$.*

(iii) *Suppose further that for all $(i, j)$, $X_{i,j}$ has mean 0 and finite moments of order $4k(\eta)$, with $k(\eta) \geq (1.5 + \varepsilon + \eta)/(2\eta)$ and $k(\eta) \in \mathbb{N}$, for some $\varepsilon > 0$. Assume that $k(\eta) \geq (2 + \varepsilon + \beta)/(2\delta_0)$.*

*Call*

$$S_p = \frac{1}{n} \sum_{i=1}^n X_i X_i'.$$

*Call $T_\alpha(S_p)$ the matrix obtained from thresholding the entries of $S_p$ at the level $Kn^{-\alpha}$ with $\alpha = 1/2 - \delta > \alpha_0$, $\delta > 0$ and $K > 0$. Then we have, if we call $\Delta_p = T_\alpha(S_p) - \Sigma_p$,*

$$\|\Delta_p\|_2 = \|T_\alpha(S_p) - \Sigma_p\|_2 \to 0 \qquad a.s. \ as \ n \to \infty,$$

*where $\|M\|_2$ is the spectral norm of the matrix $M$.*

We postpone a short discussion of the meaning of this theorem to after the statement of Theorem 2, which is arguably more interesting practically.

PROOF OF THEOREM 1. We divide the proof into two parts. The first part consists in showing the "oracle" version of the theorem, that is, showing that operator norm consistency happens when one is given the pairs $(i, j)$ for which $\sigma_p(i, j) = 0$. The second part shows that empirical thresholding does not affect this result.

Let us first remind the reader of a variant of Hölder's inequality. Let $A_1, \ldots, A_m$ be random variables with finite absolute $m$th moment. Then we have

$$\left| \mathbf{E} \left( \prod_{i=1}^{m} A_i \right) \right| \leq \prod_{i=1}^{m} \mathbf{E}(|A_i|^m)^{1/m}.$$

Note that for the case $m = 2$, this is just the Cauchy–Schwarz inequality. So the result is true when $m = 2$. We prove it by induction on $m$. Suppose therefore it is true for all integers less than or equal to $m - 1$. Call $B_1 = \prod_{i=2}^{m} A_i$. By Hölder's inequality, we have

$$|\mathbf{E}(A_1 B_1)| \leq (\mathbf{E}(|A_1|^m))^{1/m} [\mathbf{E}(|B_1|^{m/(m-1)})]^{(m-1)/m}.$$

Now, by the induction hypothesis, applied to the random variables $|A_i|^{m/(m-1)}$,

$$\mathbf{E}(|B_1|^{m/(m-1)}) = \mathbf{E} \left( \prod_{i=2}^{m} |A_i|^{m/(m-1)} \right) \leq \prod_{i=2}^{m} [\mathbf{E}(|A_i|^m)]^{1/(m-1)}.$$

Therefore, $[\mathbf{E}(|B_1|^{m/(m-1)})]^{(m-1)/m} \leq \prod_{i=2}^{m} \mathbf{E}(|A_i|^m)^{1/m}$ and the inequality is verified.

Now given $\gamma(2k)$, a closed walk of length $2k$ and the associated matrix $M$, we clearly have

$$(1) \qquad |\mathbf{E}(w_{\gamma(2k)})| \leq \mathbf{E}(|w_{\gamma(2k)}|) \leq \prod_{j=1}^{2k} [\mathbf{E}(|M(i_j, i_{j+1})|^{2k})]^{1/(2k)},$$

assuming for a moment that the relevant moments exist.

**Oracle part of the proof.** Let us first introduce some notation. We denote by $\sigma_p(i, j)$ the $(i, j)$th entry of $\Sigma_p$, the population covariance. We call oracle$(S_p)$ the matrix with entries $S_p(i, j) \mathbf{1}_{\sigma_p(i, j) \neq 0}$ and

$$\Xi_p = \text{oracle}(S_p) - \Sigma_p.$$

Note that we have

$$\Xi_p(i, j) = (S_p(i, j) - \sigma_p(i, j)) \mathbf{1}_{\sigma_p(i, j) \neq 0}.$$

In the oracle setting, where we assume we know the patterns of zeros in $\Sigma_p$, so we focus on the matrix $\Xi_p$. Clearly $\Sigma_p$ and $\Xi_p$ have the same patterns of 0's and nonzero, and so if $\Sigma_p$ is $\beta$-sparse, so is $\Xi_p$. Equation (1) shows that if we can control the moments $(\Xi_p(i,j))^{2k}$, we will be able to bound the expected weight of each walk. Now we remark that we can write

$$\Xi_p(i,j) = \frac{1}{n}\sum_{m=1}^{n} Z_m,$$

where $Z_m$'s are independent, identically distributed and with mean 0, since $S_p$ is unbiased for $\Sigma_p$. By expanding the power, we get that

$$(\Xi_p(i,j))^{2k} = \frac{1}{n^{2k}}\sum_{i_1,\ldots,i_{2k}} Z_{i_1}\cdots Z_{i_{2k}}.$$

This last quantity can be rewritten

$$Z_{i_1}\cdots Z_{i_{2k}} = \prod_{i=1}^{n} Z_i^{k_i} \qquad \text{with } \sum_{i=1}^{n} k_i = 2k \quad \text{and} \quad k_i \geq 0.$$

We now remark that if there exists $i_0$ such that $k_{i_0}=1$, then $\mathbf{E}(Z_{i_1}\cdots Z_{i_{2k}})=0$, by independence and the fact that each of the $Z_i$'s have mean 0. Therefore, in the expansion of $(\Xi_p(i,j))^{2k}$, only the terms for which all nonzero $k_i$'s are greater than or equal to 2 will contribute to the expectation. Counting the number of distinct indices appearing in the product above allows us to get a first-order estimate of $\mathbf{E}(\Xi_p(i,j)^{2k})$. As a matter of fact, the contribution of products with $q$ distinct indices is of order $n^{-2k}n^q$, by simply counting how many such products there are. So we see that to first order, the only products that matter are those for which all the $Z_i$'s raised to a nonzero power are raised to the power 2. Denoting by $n^{[k]}$ the $k$th factorial moment $n(n-1)\cdots(n-k+1)$, we have, assuming that $\mathbf{E}(Z_i^{2k}) < \infty$,

$$\mathbf{E}(\Xi_p(i,j))^{2k} \leq \frac{n^{[k]}}{n^{2k}}\frac{2k!}{2^k k!}[\mathbf{E}(Z_i^2)]^k + \frac{1}{n^{2k}}\mathrm{O}(n^{k-1}) = \mathrm{O}\left(\frac{1}{n^k}\right),$$

$$\text{if } k \text{ is fixed and } n \to \infty.$$

We therefore have

$$[\mathbf{E}(\Xi_p(i,j))^{2k}]^{1/(2k)} = \mathrm{O}\left(\frac{1}{\sqrt{n}}\right).$$

In particular, the weight of a closed walk of length $2k$ on the graph with adjacency matrix $A_p(\Sigma_p)$ [or $A_p(\Xi_p)$] and weights $\Xi_p(i,j)$ has the property that

$$|\mathbf{E}(w_{\gamma(2k)})| \leq \mathbf{E}(|w_{\gamma(2k)}|) = \mathrm{O}(n^{-k}).$$

Since we have assumed that $\Sigma_p$ and therefore $\Xi_p$ are $\beta$-sparse, we have

$$\mathbf{E}(\mathrm{trace}(\Xi_p^{2k})) = \mathrm{O}(p^{1+\beta(2k-1)}n^{-k}).$$

Since our assumption $p \asymp n$ implies that $p/n$ remains bounded, we see that $\mathbf{E}(\mathrm{trace}(\Xi_p^{2k})) = \mathrm{O}(n^{1/2+\eta-2k\eta})$, where $\eta = 1/2 - \beta$. In particular, if $k$ is chosen such that

$$k \geq \frac{1.5 + \varepsilon + \eta}{2\eta},$$

we see that

$$\mathbf{E}(\|\!|\Xi_p|\!\|_2^{2k}) \leq \mathbf{E}(\mathrm{trace}(\Xi_p^{2k})) = \mathrm{O}(n^{-(1+\varepsilon)}),$$

because $\Xi_p$ is a symmetric matrix, so its spectral norm squared is one of its eigenvalues squared. Using Chebyshev's inequality and the first Borel–Cantelli lemma, we conclude that

$$\|\!|\Xi_p|\!\|_2 \to 0 \qquad \text{a.s.}$$

Note that $2k > 1 + 1/(2\eta)$ would have guaranteed convergence in probability. The above proof is correct if $Z_m$ has a finite $2k$th moment. Since $Z_m = X_{m,i}X_{m,j}$, the assumption that the entries of the data matrix $X$ have a $4k$th moment guarantees the existence of a $2k$th moment for $Z_m$, using, for instance, the Cauchy–Schwarz inequality.

We have shown that $\|\!|\mathrm{oracle}(S_p) - \Sigma_p|\!\|_2 \to 0$ almost surely, when the conditions of the theorem are satisfied.

**Nonoracle part of the proof.** We now turn to the nonoracle version of the procedure. It is clear that all we need to do at this point is to show that the thresholding procedure will lead a.s. to the right adjacency matrix. Recall the notation $\Delta_p = T_\alpha(S_p) - \Sigma_p$, the difference between our estimator and the population covariance. Call $B_p$ the event $B_p = \{$at least one mistake is made by thresholding$\}$, that is, $A_p(T_\alpha(S_p)) \neq A_p(\Sigma_p)$. Call $E_p$ the event $\{\|\!|\Delta_p|\!\|_2 > \varepsilon\}$ and $F_p$ the event $\{\|\!|\Xi_p|\!\|_2 > \varepsilon\}$ (we do not index these events by $\varepsilon$ to alleviate the notation). Note that

$$E_p = (E_p \cap B_p) \cup (E_p \cap B_p^c) \subseteq B_p \cup (E_p \cap B_p^c) = B_p \cup (F_p \cap B_p^c) \subseteq B_p \cup F_p.$$

We have already seen that $P(F_p \text{ infinitely often}) = 0$, so if we can show that $P(B_p \text{ i.o.}) = 0$, we will have $P(E_p \text{ i.o.}) = 0$, as desired.

Call $O_p = \mathrm{oracle}(S_p)$ and $S_p^- = S_p - O_p$, where $\mathrm{oracle}(S_p)$ is defined above. Note that $S_p^-$ has nonzero entries only where $\Sigma_p$ has entries equal to 0; when that is the case, $S_p^-(i,j) = \hat{\sigma}(i,j)$. We call $\mathcal{D}_p$ the set of pairs $(i,j)$ such that $\sigma(i,j) = 0$, that is,

$$\mathcal{D}_p = \{(i,j): \sigma_p(i,j) = 0\}.$$

We will first show that the maximal element of $S_p^-$ stays below $n^{-\alpha}$ a.s. Note that in general, for a random matrix $M$ and index sets $I$ and $J$,

$$P\Big(\max_{i\in I, j\in J}|m_{i,j}| > \varepsilon\Big) \leq \sum_{i\in I, j\in J} P(|m_{i,j}| > \varepsilon).$$

The same moment computations as the ones we made for $\Xi_p$ above show that for the elements of $S_p$ corresponding to $\sigma_p(i,j) = 0$, we have $\mathbf{E}(S_p(i,j)^{2k}) = O(n^{-k})$. Therefore,

$$P\Big(\max_{\mathcal{D}_p}|S_p(i,j)| > Cn^{-\alpha}\Big) \leq \sum_{(i,j)\in\mathcal{D}_p} \mathbf{E}(S_p(i,j)^{2k})\frac{n^{2k\alpha}}{C^{2k}} = O(p^2 n^{2k\alpha} n^{-k}).$$

Since we assumed that $p \asymp n$, we see that if $k(1-2\alpha) - 2 \geq 1 + \varepsilon$,

$$P\Big(\max_{\mathcal{D}_p}|S_p(i,j)| > Cn^{-\alpha} \text{ i.o.}\Big) = 0,$$

by the first Borel–Cantelli lemma. In other words, if we call $(T_\alpha(S_p))^-$ the thresholded version of the part of $S_p$ that corresponds to indices in $\mathcal{D}_p$, we have that $P((T_\alpha(S_p))^- \neq 0 \text{ i.o.}) = 0$.

We now turn our attention to $\mathcal{D}_p^c$, that is, the set of indices for which $\sigma_p(i,j) \neq 0$. Recall that we assumed that these $\sigma_p(i,j)$ satisfied $|\sigma_p(i,j)| \geq Cn^{-\alpha_0}$ and $\alpha_0 < \alpha$. Now note for $(i,j)$ in $\mathcal{D}_p^c$, and $\sigma_p(i,j) \geq 0$, we have $\{|S_p(i,j)| < Cn^{-\alpha}\} \subseteq \{0 \leq \sigma_p(i,j) - Cn^{-\alpha} \leq \sigma_p(i,j) - S_p(i,j)\}$. So, in this case, by using the moment computations made above, and using $C$ to denote a generic constant, we have

$$P(|S_p(i,j)| < Cn^{-\alpha}) \leq \frac{\mathbf{E}(\sigma_p(i,j) - S_p(i,j))^{2k}}{(\sigma_p(i,j) - Cn^{-\alpha})^{2k}} = O(n^{-k} n^{2k\alpha_0}).$$

Similarly, when $\sigma_p(i,j) < 0$, we have

$$P(|S_p(i,j)| < Cn^{-\alpha}) \leq \frac{\mathbf{E}(\sigma_p(i,j) - S_p(i,j))^{2k}}{(|\sigma_p(i,j)| - Cn^{-\alpha})^{2k}} = O(n^{-k} n^{2k\alpha_0}).$$

Now note that because $\Sigma_p$ is $\beta$-sparse, there are at most $O(p^{1+\beta})$ nonzero coefficients in $\Sigma_p$; indeed $\phi_p(2)$ counts the number of nonzero coefficients in $\Sigma_p$. From this we conclude that

$$P(\exists (i_0,j_0)\in\mathcal{D}_p^c : |S_p(i_0,j_0)| < Cn^{-\alpha}) \leq O(n^{k(2\alpha_0-1)} p^{1+\beta}).$$

So if

$$k \geq \frac{2+\varepsilon+\beta}{1-2\alpha_0} = \frac{2+\varepsilon+\beta}{2\delta_0},$$

then, almost surely, no $S_p(i,j)$ will be wrongly thresholded, if $(i,j) \in \mathcal{D}_p^c$. Combining this result with the one on the indices in $\mathcal{D}_p$, we have

$$P(B_p \text{ i.o.}) = 0,$$

and we have the result announced in the theorem. $\square$

It is, however, more common practice to use as our estimator of the covariance matrix the sample covariance matrix that differs slightly from the matrix $S_p$ used above, which is the maximum likelihood estimator in the (mean 0) Gaussian case. We now show that for the usual estimator the same strategy works.

THEOREM 2 (Sample covariance matrix).  *Suppose the assumptions of Theorem 1 are satisfied, but allow now $X_i$ to have a nonzero mean $\mu$. Call*

$$S_p = \frac{1}{n-1} \sum_{i=1}^{n} (X_i - \bar{X})(X_i - \bar{X})'.$$

*Then the result of Theorem 1 holds; namely, the matrix $T_\alpha(S_p) - \Sigma_p$ converges a.s. in spectral norm to 0.*

The previous theorem basically means that if the covariance matrix $\Sigma_p$ is sparse enough, and if the data come from a distribution with enough moments, then thresholding the sample covariance matrix by keeping only terms that are a bit larger than $1/\sqrt{n}$ is a good idea and will lead to an estimator that is consistent in operator norm. This is in stark contrast to simply using the sample covariance matrix, when in the asymptotics considered here, we would not, in general (e.g., as soon as $p/n \to l \neq 0$), have consistency even at the level of the vector of eigenvalues; in the case of $\Sigma_p = \text{Id}$, this is a consequence of the results of [13] or [21] and we refer to [12] for a thorough discussion.

PROOF OF THEOREM 2.  The proof proceeds as the one of Theorem 1. Since $S_p$ is still unbiased for $\Sigma_p$, the only thing we have to show here is that the $2k$th central moments of $S_p(i,j)$ decay in the same fashion as they did in Theorem 1. Firstly let us note that

$$S_p(i,j) = \frac{1}{n-1} \sum_{l=1}^{n} (X_{l,i} - \mu_i)(X_{l,j} - \mu_j) - \frac{n}{n-1} (\bar{X}_i - \mu_i)(\bar{X}_j - \mu_j),$$

so

$$S_p(i,j) - \sigma_p(i,j) = \frac{1}{n-1} \sum_{l=1}^{n} ((X_{l,i} - \mu_i)(X_{l,j} - \mu_j) - \sigma_p(i,j))$$
$$- \frac{n}{n-1} \left( (\bar{X}_i - \mu_i)(\bar{X}_j - \mu_j) - \frac{1}{n}\sigma_p(i,j) \right).$$

Now, since $(a + b)^{2k} \leq 2^{2k}(a^{2k} + b^{2k})$, we see that we will have the result we need if we can bound each term in the right-hand side of the previous equation. The technique we used above immediately shows that

$$\mathbf{E}\left(\frac{1}{n-1}\sum_{l=1}^n [(X_{l,i} - \mu_i)(X_{l,j} - \mu_j) - \sigma_p(i,j)]\right)^{2k} = \mathrm{O}\left(\frac{1}{n^k}\right),$$

assuming for a moment that all the needed moments exist. For the other part of the equation, the same argument shows that the only thing we need to control is $\mathbf{E}((\bar{X}_i - \mu_i)(\bar{X}_j - \mu_j))^{2k}$, since the assumptions we made about the moments of $X_i$ guarantee that $\sigma_p(i,j)$ is bounded in $p$. Using the Cauchy–Schwarz inequality, it is clear that all we need to do is control $\mathbf{E}(\bar{X}_i - \mu_i)^{4k}$, for all $i$. But $\bar{X}_i - \mu_i$ is a sum of independent mean-0 random variables and the computations we made in the proof of Theorem 1 show that

$$\mathbf{E}(\bar{X}_i - \mu_i)^{4k} = \mathrm{O}\left(\frac{1}{n^{2k}}\right).$$

Therefore,

$$\mathbf{E}((\bar{X}_i - \mu_i)(\bar{X}_j - \mu_j))^{2k} = \mathrm{O}\left(\frac{1}{n^{2k}}\right).$$

So we conclude that

$$\mathbf{E}(S_p(i,j) - \sigma_p(i,j))^{2k} = \mathrm{O}\left(\frac{1}{n^k}\right),$$

just as in the case of the Gaussian MLE estimator. This is all we need to complete the proof of Theorem 2, since the last steps follow exactly from the proof of Theorem 1. The assumptions made guarantee that all the moments used above exist and are finite. □

We note that the distribution of the entries of $X$ can change with $n$ and $p$ as long as the moment conditions are satisfied and the bounds on the moments are uniform in $n$ and $p$. We now turn to the question of estimating correlation matrices.

THEOREM 3 (Correlation matrices). *Under the assumptions of Theorem 1, but requiring the boundedness of the $8k(\eta)$th moments of the $X_{i,j}$'s in assumption* (iii), *if $\Sigma_p$ is now the correlation matrix of the vector $X_i$, and if $S_p$ is now the sample correlation matrix, we have as before*

$$\|T_\alpha(S_p) - \Sigma_p\|_2 \to 0 \qquad a.s.$$

REMARK 1. We note that the moment assumption can be relaxed to $4k(\eta)$, if, for instance, one assumes that $\|\Sigma_p\|_2$ is bounded. This is a simple

consequence of the fact that, if we call $D_p$ the diagonal of the sample covariance matrix $\tilde{S}_p$, the sample correlation matrix $S_p$ is equal to $D_p^{-1/2} \tilde{S}_p D_p^{-1/2}$. Because $\Sigma_p$ has only 1-s on the diagonal, our results on operator-norm-consistent estimation of $\Sigma_p$ imply in particular that $\|\|D_p - \mathrm{Id}_p\|\|_2$ tends to 0 a.s. Elementary algebra then shows that if $\|\|\Sigma_p\|\|_2$ is, for instance, bounded, $\|\|D_p^{-1/2} \tilde{S}_p D_p^{-1/2} - \Sigma_p\|\|_2$ also tends to 0 a.s., because $\|\|\tilde{S}_p - \Sigma_p\|\|_2$ does, according to Theorem 2.

PROOF OF THEOREM 3. Because of invariance of the problem by centering and scaling, we can assume that the row vector $X_i$ has mean 0, and that the diagonal of its covariance matrix $\Sigma_p$ is full of 1. Then we have $\rho(i,j) = \sigma_p(i,j)$. From the proof of Theorem 1, it is clear that if we can show that $\mathbf{E}(S_p(i,j) - \rho(i,j))^{2k} = \mathrm{O}(n^{-k})$ for all $(i,j)$, the same technique as above will lead to the theorem. To show that this is indeed the case, we first make the following elementary remark, which prepares the study of $\hat{\rho}(i,j) - \rho(i,j)$. Suppose that $F_n$ and $G_n$ are random variables, with $\mathbf{E}(F_n - \rho)^{2k} = \mathrm{O}(n^{-k})$ for some $\rho \in [-1,1]$, $\mathbf{E}(G_n - 1)^{2k} = \mathrm{O}(n^{-k})$ and further $|F_n/G_n| \leq 1$. Call $\Omega_n(\varepsilon)$ the event $\{\omega : |G_n - 1| < \varepsilon\}$. We have

$$
\begin{aligned}
\mathbf{E}\left(\frac{F_n}{G_n} - \rho\right)^{2k} &= \mathbf{E}\left(\left[\frac{F_n}{G_n} - \rho\right]^{2k}[1_{\Omega_n(\varepsilon)} + 1_{\Omega_n^c(\varepsilon)}]\right) \\
&\leq \mathbf{E}\left(\left[\frac{F_n}{G_n} - \rho\right]^{2k} 1_{\Omega_n(\varepsilon)}\right) + \mathbf{E}\left(\left[\frac{F_n}{G_n} - \rho\right]^{2k} 1_{\Omega_n^c(\varepsilon)}\right) \\
&\leq \mathbf{E}\left(\left[\frac{F_n - \rho}{G_n} - \rho\left(1 - \frac{1}{G_n}\right)\right]^{2k} 1_{\Omega_n(\varepsilon)}\right) + 2^{2k}\mathbf{E}(1_{\Omega_n^c(\varepsilon)}) \\
&\leq 2^{2k}\mathbf{E}\left(\left[\frac{F_n - \rho}{G_n}\right]^{2k} 1_{\Omega_n(\varepsilon)} + \rho^{2k}\left[1 - \frac{1}{G_n}\right]^{2k} 1_{\Omega_n(\varepsilon)}\right) \\
&\quad + 2^{2k}\mathbf{E}(1_{\Omega_n^c(\varepsilon)}) \\
&\leq \frac{2^{2k}}{(1-\varepsilon)^{2k}}\{\mathbf{E}((F_n - \rho)^{2k} 1_{\Omega_n(\varepsilon)}) + \rho^{2k}\mathbf{E}((G_n - 1)^{2k} 1_{\Omega_n(\varepsilon)})\} \\
&\quad + 2^{2k}\mathbf{E}(1_{\Omega_n^c(\varepsilon)}).
\end{aligned}
$$

By Chebyshev's inequality, and our assumptions, it is clear that

$$
\mathbf{E}\left(\frac{F_n}{G_n} - \rho\right)^{2k} = \mathrm{O}\left(\frac{1}{n^k}\right).
$$

Now we claim that this remark applies in the case of the correlation matrix. We have

$$
\hat{\rho}(i,j) = \frac{F_n(i,j)}{G_n(i,j)},
$$

where

$$F_n(i,j) = \frac{1}{n-1} \sum_{l=1}^{n} (X_{l,i} - \bar{X}_i)(X_{l,j} - \bar{X}_j),$$

and $G_n(i,j) = \sqrt{F_n(i,i)F_n(j,j)}$.

From the moment computations made in the proof of Theorem 2, we see that we have $\mathbf{E}(F_n(i,j) - \rho(i,j))^{2k} = \mathrm{O}(n^{-k})$, for all $(i,j)$. Let us denote by $Y_n(i) = F_n(i,i)$; this result implies that

$$\mathbf{E}(\sqrt{Y_n(i)} - \sqrt{\rho(i,i)})^{2k} \leq \mathbf{E}(Y_n(i) - \rho(i,i))^{2k}/\rho(i,i)^k = \mathbf{E}(Y_n(i)-1)^{2k}$$
$$= \mathrm{O}(n^{-k}).$$

If we denote $\alpha_n = \sqrt{Y_n(i)}$ and $\beta_n = \sqrt{Y_n(j)}$, we have, since $\alpha_n \beta_n - 1 = (\alpha_n - 1)\beta_n + \beta_n - 1$,

$$\mathbf{E}(\alpha_n \beta_n - 1)^{2k} \leq 2^{2k}[\mathbf{E}(\beta_n - 1)^{2k} + \mathbf{E}(\beta_n^{2k}(\alpha_n - 1)^{2k})]$$
$$\leq 2^{2k}[\mathbf{E}(\beta_n - 1)^{2k} + \sqrt{\mathbf{E}(\beta_n^{4k})}\sqrt{\mathbf{E}(\alpha_n - 1)^{4k}}].$$

We have already seen that $\mathbf{E}(\beta_n - 1)^{2k} = \mathrm{O}(n^{-k})$, and since we are assuming the existence of a $8k$th moment for the $X_{l,j}$, we also have $\sqrt{\mathbf{E}(\alpha_n - 1)^{4k}} = \mathrm{O}(n^{-k})$. To conclude that $\mathbf{E}(\alpha_n \beta_n - 1)^{2k} = \mathrm{O}(n^{-k})$, we just need to show that $\mathbf{E}(\beta_n)^{4k}$ is bounded. But $\beta_n^{4k} = (\beta_n - 1 + 1)^{4k} \leq 2^{4k}((\beta_n - 1)^{4k} + 1)$, from which we conclude that $\mathbf{E}(\beta_n)^{4k}$ is bounded, since $(\beta_n - 1)^{4k} = \mathrm{O}(n^{-2k})$. We now see that $G_n(i,j) = \sqrt{Y_n(i)Y_n(j)}$ satisfies with $F_n$ the conditions needed to conclude that

$$\mathbf{E}\left(\frac{F_n(i,j)}{G_n(i,j)} - \rho(i,j)\right)^{2k} = \mathbf{E}(\widehat{\rho}(i,j) - \rho(i,j))^{2k} = \mathrm{O}(n^{-k}). \qquad \square$$

3.2. *Approximation of nonsparse matrices by sparse matrices.* It is natural to ask whether a thresholding approach can also lead to good results when dealing with matrices which are not sparse per se, but have many coefficients close to zero. In other words, we would like to know when we can approximate nonsparse matrices by sparse matrices obtained by thresholding a sample covariance or correlation matrix. We now present two propositions that relax the sparsity assumptions and still lead to spectral norm convergence. The most general one basically says that if the population covariance matrix can be approximated by a sparse matrix and does not have too many coefficients close to the threshold level $1/\sqrt{n}$, then estimating the (not necessarily sparse) population covariance by thresholding the sample covariance matrix will lead to good results.

Here is our first result in this direction:

PROPOSITION 1. *Making the same general assumptions as in the theorems above [i.e., assumptions* (i)–(iii) *in Theorem 1 are excluded], we now assume that:*

- *There exists $T_{\alpha_1}(\Sigma_p) = \widetilde{\Sigma}_p$, a version of $\Sigma_p$ thresholded at $Cn^{-\alpha_1}$, that is, $\beta$-sparse, with $\beta = 1/2 - \eta$ and $\eta > 0$. Further we assume that $\|\widetilde{\Sigma}_p - \Sigma_p\|_2 \to 0$.*
- *We call $I_{\alpha_1, \alpha_0}$ the set of indices of those $\sigma(i, j)$ for which $Cn^{-\alpha_1} < |\sigma(i, j)| < Cn^{-\alpha_0}$, with $\alpha_0 < \alpha_1 \leq 1/2 - \delta_0$, for some $\delta_0 > 0$.*
- *The adjacency matrix corresponding to $I_{\alpha_1, \alpha_0}$ is $\gamma$-sparse, for some $\gamma \leq \alpha_0 - \zeta_0$, where $\zeta_0 > 0$.*
- *The random variables $X_{i,j}$ have moments of order $4k$ ($8k$ in the correlation case), with $k$ satisfying the assumptions put forth in Theorem 1, assumption* (iii), *as well as $k \geq (2 + \varepsilon - \gamma)/(1 - 2\gamma)$, for some $\varepsilon > 0$.*

*Now if we choose $\alpha \in (\alpha_0, \alpha_1)$, then the conclusions of all the theorems above apply:*

$$\|T_\alpha(S_p) - \Sigma_p\|_2 \to 0 \qquad a.s., \ as \ n \to \infty.$$

While this proposition might appear full of hard-to-check assumptions, we believe it is useful and not so hard to use when checking whether thresholding is a reasonable idea for a particular estimation problem. We give an example after stating Fact 3 below. Finally, we note that under the assumptions stated, both $T_{\alpha_0}(\Sigma_p)$ and $T_{\alpha_1}(\Sigma_p)$ are good approximations of $\Sigma_p$ in operator norm.

PROOF OF PROPOSITION 1. In the proof we assume without loss of generality that $C = 1$, which allows us to avoid cumbersome notation. [As the reader will see, replacing $n^{-\alpha_0}$ by $Cn^{-\alpha_0}$ every time an $n^{-\alpha_0}$ (and similarly for $n^{-\alpha_1}$) appears does not change anything in the proof.]

From the previous proofs, we see that we can divide

$$T_\alpha(S_p) = M_0 + M_1 + M_2$$

into three parts. $M_0$ corresponds to the indices $(i, j)$ for which $\sigma(i, j)$ is larger (in absolute value) than $n^{-\alpha_0}$, $M_1$ corresponds to indices in $I_{\alpha_0, \alpha_1}$, and $M_2$ to those indices for which $|\sigma(i, j)| < n^{-\alpha_1}$. Similarly, we can write with the same partition of indices,

$$\Sigma_p = T_{\alpha_0}(\Sigma_p) + [T_{\alpha_1}(\Sigma_p) - T_{\alpha_0}(\Sigma_p)] + [\Sigma_p - T_{\alpha_1}(\Sigma_p)] = \Sigma_0 + \Sigma_1 + \Sigma_2.$$

With the same notation for the subparts of $\Sigma$, we have from the computations we made in the proofs of the previous theorems that $\|M_0 - \Sigma_0\|_2 \to 0$ a.s. (by the oracle part of the proofs), and $\|M_2\|_2 \to 0$ a.s., since the $\widehat{\sigma}_p(i, j)$ corresponding to $|\sigma(i,j)| < n^{-\alpha_1}$ will all be (a.s.) thresholded to 0 if the thresholding level is $n^{-\alpha}$, $\alpha < \alpha_1$.

Note that $\Sigma_0 + \Sigma_1 = \widetilde{\Sigma}_p$, so $\|\|\Sigma_2\|\|_2 \to 0$. To reach the conclusions of the proposition, we need to show that we control $M_1 - \Sigma_1$ in operator norm.

Recall that our assumption is that $\Sigma_1$ is $\gamma$-sparse. We call

$$\Sigma_1 = T_\alpha(\Sigma_1) + R_\alpha(\Sigma_1),$$

where $T_\alpha(\Sigma_1)$ is the version of $\Sigma_1$ thresholded at $n^{-\alpha}$. It is of course also $\gamma$-sparse and so is $R_\alpha(\Sigma_1)$. This implies that

$$\|\|R_\alpha(\Sigma_1)\|\|_2^{2k} \leq \mathrm{trace}((R_\alpha(\Sigma_1))^{2k}) \leq f(k)p^{\gamma(2k-1)}pn^{-2k\alpha},$$

which goes to 0 if $\gamma \leq \alpha - \varepsilon$. So we can find $k_0$, an integer independent of $n$ and $p$, such that the right-hand side goes to 0 as $n$ and $p$ go to infinity. This implies that $\|\|R_\alpha(\Sigma_1)\|\|_2 \to 0$, as $n$ tends to infinity.

Using the oracle proof of Theorem 1, we see that if we make no error in thresholding for the indices in $I_{\alpha_0,\alpha_1}$, then $\|\|\mathrm{oracle}_\alpha(M_1) - T_\alpha(\Sigma_1)\|\|_2$ tends to 0 a.s. Therefore, all we need to do is check that we control the operator norm of the matrix of possible errors, that is, the difference $M_1 - \mathrm{oracle}_\alpha(M_1)$. Let us call $\Upsilon_1$ this matrix of potential errors. There are two types of possible errors: either a coefficient is thresholded when it should not have been, or it is not thresholded when it should have been thresholded. So

$$\Upsilon_1(i,j) = \begin{cases} 0, & \text{if correctly thresholded } \hat{\sigma}(i,j), \\ \hat{\sigma}(i,j), & \text{if } |\sigma_{i,j}| \leq n^{-\alpha} \text{ but did not threshold in } M_1, \\ -\hat{\sigma}(i,j), & \text{if } |\sigma_{i,j}| > n^{-\alpha} \text{ but did threshold in } M_1. \end{cases}$$

In any case, we conclude that $|\Upsilon_1(i,j)| \leq |\hat{\sigma}(i,j)| \leq |\hat{\sigma}(i,j) - \sigma(i,j)| + |\sigma(i,j)|$. Let us call $\Upsilon_{11}$ the matrix

$$\Upsilon_{11}(i,j) = 1_{\Upsilon_1(i,j)\neq 0}|\hat{\sigma}(i,j) - \sigma(i,j)|,$$

and $\Upsilon_{12}$ the matrix with entries

$$\Upsilon_{12}(i,j) = 1_{\Upsilon_1(i,j)\neq 0}|\sigma(i,j)|.$$

Note that all the indices where $\Upsilon_1$ has potentially nonzero entries are in $I_{\alpha_0,\alpha_1}$, so the corresponding adjacency matrix is $\gamma$-sparse. Clearly, the same is true for $\Upsilon_{11}$ and $\Upsilon_{12}$.

Now, $\|\|\Upsilon_1\|\|_2 \leq \|\|\Upsilon_{11} + \Upsilon_{12}\|\|_2$, according to Lemma A.2 in the Appendix. Therefore,

$$\|\|\Upsilon_1\|\|_2 \leq \|\|\Upsilon_{11}\|\|_2 + \|\|\Upsilon_{12}\|\|_2.$$

Using the fact that all the entries of $\Upsilon_{12}$ are less than $n^{-\alpha_0}$ in absolute value, and the fact that $\Upsilon_{12}$ is $\gamma$-sparse, we have, for $k$ integer, according to Lemma A.1 below,

$$\|\|\Upsilon_{12}\|\|_2^{2k} \leq \mathrm{trace}(\Upsilon_{12}^{2k}) = \mathrm{O}(p^{\gamma(2k-1)+1}n^{-2k\alpha_0}).$$

So since $p \asymp n$, and we assumed that $\gamma < \alpha_0$, we see that we can find $k$ (finite) such that the right-hand side goes to 0 as $n \to \infty$. Actually, any $k > (1 - \gamma)/(2(\alpha_0 - \gamma))$ is a valid choice. So we conclude that

$$\|\Upsilon_{12}\|_2 \to 0 \qquad \text{as } n \to \infty.$$

On the other hand, since

$$\mathbf{E}(\Upsilon_{11}(i,j))^{2k} \leq \mathbf{E}(\hat{\sigma}(i,j) - \sigma(i,j))^{2k} = \mathrm{O}(n^{-k}),$$

we conclude as before that the expected weight of a walk ("on" $\Upsilon_{11}$) of length $2k$ is $\mathrm{O}(n^{-k})$. Using the assumption of $\gamma$-sparsity of the matrix $\Sigma_1$, we conclude that $\mathbf{E}(\mathrm{trace}(\Upsilon_{12}^{2k}))$ is $\mathrm{O}(p^{\gamma(2k-1)+1} n^{-k})$. Therefore, if we can pick $k \geq (2 + \varepsilon - \gamma)/(1 - 2\gamma)$, and finite, we have a.s. convergence of $\|\Upsilon_{11}\|_2$ to 0. Now note that our moment requirements imply that indeed we can pick $k$ (finite), with the property that $k \geq (2 + \varepsilon - \gamma)/(1 - 2\gamma)$. Hence,

$$\|\Upsilon_1\|_2 \to 0 \qquad \text{a.s.}$$

This concludes the proof since we have bounded $\|T_\alpha(S_p) - \Sigma_p\|_2$ by a sum of operator norms of matrices, all of which are going to 0 a.s. $\quad\square$

We also have the following proposition.

PROPOSITION 2. *Let us denote by $|\Sigma_p|_{\mathrm{Had}}$ the Hadamard absolute value of $\Sigma_p$, that is, the matrix whose $(i,j)$th entry is equal to $|\sigma(i,j)|$. Suppose that $\||\Sigma_p|_{\mathrm{Had}}\|_2$ is uniformly bounded in $p$. Suppose $p \asymp n$. Let us call, for $\alpha > 0$,*

$$R_\alpha(\Sigma_p) = \Sigma_p - T_\alpha(\Sigma_p).$$

*Suppose that $\||R_{\alpha_0}(\Sigma_p)|_{\mathrm{Had}}\|_2 \to 0$, for some given $\alpha_0$, with $\alpha_0 < 1/2 - \delta_0$, and $\delta_0 > 0$. Then, for $\alpha_0 < \alpha < 1/2 - \delta_0$, we have*

$$\|T_\alpha(S_p) - \Sigma_p\|_2 \to 0 \qquad a.s.,$$

*provided the moment conditions in Theorem 1(iii) are satisfied, with the parameters $\delta_0$, $\eta = \delta_0$ and hence $\beta = 1/2 - \delta_0$.*

PROOF. Since we have assumed that $\||\Sigma_p|_{\mathrm{Had}}\|_2 < \infty$, we therefore have $\||T_\alpha(\Sigma_p)|_{\mathrm{Had}}\|_2 < \infty$, by using Lemma A.2 in the Appendix. Now since the smallest nonzero entry of $|T_\alpha(\Sigma_p)|_{\mathrm{Had}}$ is greater than $n^{-\alpha}$, we also have, by the same arguments as those developed in Lemma A.2, for any integer $k$,

$$\mathrm{trace}((|T_\alpha(\Sigma_p)|_{\mathrm{Had}})^k) \geq \phi_p(k) n^{-\alpha k}.$$

Hence,

$$\mathrm{trace}((|T_\alpha(\Sigma_p)|_{\mathrm{Had}})^{2k})^{1/(2k)} \geq (\phi_p(2k))^{1/(2k)} n^{-\alpha}.$$

We now note that without loss of generality, we can assume in our definition of $\beta$-sparsity that $f(k) \geq 1$. Hence, under our assumptions, $T_\alpha(\Sigma_p)$ has to be at most $\alpha$-sparse, for otherwise the lower bound in the previous equation would go to infinity, and we just saw that $|||T_\alpha(\Sigma_p)|_{\mathrm{Had}}||_2$ is bounded. [Recall that $|||T_\alpha(\Sigma_p)|_{\mathrm{Had}}||_2 = \lim_{k\to\infty} \mathrm{trace}((|T_\alpha(\Sigma_p)|_{\mathrm{Had}})^{2k})^{1/(2k)}$.]

Now, it is clear that if $\alpha_1 > \alpha_0$, and therefore $n^{-\alpha_1} < n^{-\alpha_0}$, $|R_{\alpha_1}(i,j)| \leq |R_{\alpha_0}(i,j)|$. This allows us to conclude that

$$|||R_{\alpha_1}|_{\mathrm{Had}}||_2 \leq |||R_{\alpha_0}|_{\mathrm{Had}}||_2.$$

Since we have assumed that $|||R_{\alpha_0}|_{\mathrm{Had}}||_2$ tends to 0, it is also the case for $|||R_\alpha|_{\mathrm{Had}}||_2$ for $\alpha \geq \alpha_0$. Since we assumed that $\alpha_0 < 1/2 - \delta_0$, we can find $\alpha_1$ such that $\alpha_0 < \alpha_1 \leq 1/2 - \delta_0$. Let us pick one such $\alpha_1$. Let us also pick $\alpha$ in $(\alpha_0, \alpha_1)$. The situation is now fairly similar to that of Proposition 1, but we cannot immediately apply this result because our control of the sparsity of $I_{\alpha_1,\alpha_0}$ is not very good.

However, we can apply similar arguments that we detail here. We use the same decompositions and notation as in the proof of this theorem. Note that $\Sigma_1 = T_{\alpha_1}(\Sigma_p) - T_{\alpha_0}(\Sigma_p)$ is a submatrix of $\Sigma$. Hence, $|R_\alpha(\Sigma_1)(i,j)| \leq |R_\alpha(\Sigma)(i,j)|$. Since $|||R_\alpha(\Sigma)|_{\mathrm{Had}}||_2$ goes to 0, we have

$$|||R_\alpha(\Sigma_1)|_{\mathrm{Had}}||_2 \to 0.$$

Note also that $||\Sigma_2||_2 = ||R_{\alpha_1}(\Sigma)||_2 \leq |||R_{\alpha_1}(\Sigma)|_{\mathrm{Had}}||_2 \to 0$. So all we need to do to complete the proof is to control the matrix $\Upsilon_1$ of potential errors. Recall that

$$|\Upsilon_1(i,j)| \leq |\hat{\sigma}(i,j)| \leq |\hat{\sigma}(i,j) - \sigma(i,j)| + |\sigma(i,j)|.$$

Recall also that all the indices where $\Upsilon_1$ has potentially nonzero entries correspond to the entries of $\Sigma_p$ whose absolute values are between $n^{-\alpha_1}$ and $n^{-\alpha_0}$, so $\Upsilon_1$ is at most $\alpha_1$-sparse. Let us call, as before, $\Upsilon_{11}$ the matrix

$$\Upsilon_{11}(i,j) = 1_{\Upsilon_1(i,j)\neq 0}|\hat{\sigma}(i,j) - \sigma(i,j)|,$$

and $\Upsilon_{12}$ the matrix with entries

$$\Upsilon_{12}(i,j) = 1_{\Upsilon_1(i,j)\neq 0}|\sigma(i,j)|.$$

Clearly,

$$||\Upsilon_1||_2 \leq |||\Upsilon_1|_{\mathrm{Had}}||_2 \leq ||\Upsilon_{11} + \Upsilon_{12}||_2 \leq |||\Upsilon_{11}||_2 + ||\Upsilon_{12}||_2.$$

Because $\Upsilon_{11}$ is $\alpha_1$-sparse and $\alpha_1 \leq 1/2 - \delta_0$, and because of the oracle part of the proof of Theorem 1, we see that, because of our moment assumptions,

$$||\Upsilon_{11}||_2 \to 0 \qquad \text{a.s.}$$

On the other hand, $\Upsilon_{12}$ is a submatrix of $|\Sigma_1|_{\mathrm{Had}}$, so

$$\| \Upsilon_{12} \|_2 \leq \| |T_{\alpha_1}(\Sigma_p) - T_{\alpha_0}(\Sigma_p)|_{\mathrm{Had}} \|_2$$
$$\leq \| |R_{\alpha_1}(\Sigma_p) - R_{\alpha_0}(\Sigma_p)|_{\mathrm{Had}} \|_2$$
$$\leq \| |R_{\alpha_1}(\Sigma_p)|_{\mathrm{Had}} \|_2 + \| |R_{\alpha_0}(\Sigma_p)|_{\mathrm{Had}} \|_2.$$

We conclude that $\| \Upsilon_{12} \|_2 \to 0$, and therefore that $\| \Upsilon_1 \|_2 \to 0$ a.s. Arguing as in the proof of Proposition 1, we finally have the result announced in Proposition 2. $\square$

The following simple fact is a clear case of applicability of the ideas of Proposition 1.

FACT 3. *Let $\varepsilon > 0$ and suppose that $T_{1+\varepsilon}(\Sigma_p)$ is $\beta$-sparse and its nonzero entries are larger in absolute value than $n^{-\alpha_0}$. Then under the same assumptions as Theorems 1, 2 and 3, we have, for $\alpha_0 < \alpha < 1/2 - \delta$,*

$$\| T_\alpha(S_p) - \Sigma_p \|_2 \to 0 \qquad a.s.$$

PROOF. Take $\alpha_1 = \alpha_0 + \delta$, where $\delta$ is small. In particular, of course, $\alpha_1 < 1/2 < 1 + \varepsilon$. Here $I_{\alpha_1, \alpha_0}$ is empty so the corresponding matrix is 0-sparse. In particular, that means that in the notation of the proof of Proposition 1, $M_1 = 0$ and similarly for $\Sigma_1$. So the results on $M_0 - \Sigma_0$ and $M_2$ apply directly and the only thing we have to check is that $\| \Sigma_2 \|_2 \to 0$. Note that $\Sigma_2$ contains only entries of order $n^{-(1+\varepsilon)}$ or smaller. Using a Frobenius norm bound, we therefore have

$$\| \Sigma_2 \|_2^2 \leq p^2 n^{-(2+2\varepsilon)} \to 0,$$

so the result is established. $\square$

*Example: a simple (permuted) Toeplitz matrix.* We consider a matrix that is often used as an example for estimation: the (Toeplitz) covariance matrix $\Sigma_p$, with $\Sigma(i,j) = \rho^{|i-j|}$, $|\rho| < 1$. Of course, we can also consider the same matrix where the variables have been randomly permuted and hence the Toeplitz structure destroyed. However, on any given line, the entries are still a (possibly random) permutation of the $\rho^{|i-j|}$. We apply Proposition 1. To do so, we just need to count how many coefficients on each row are between $n^{-\alpha_1}$ and $n^{-\alpha_0}$, for $\alpha_0$ and $\alpha_1$ to be chosen later. Note that $|\rho|^k \leq n^{-\alpha_1}$ is equivalent to $k \geq \log(n)\alpha_1 / \log(1/|\rho|)$. So $T_{\alpha_1}(\Sigma_p)$ is asymptotically 0-sparse, as it contains only $O(\log(n))$ nonzero terms on each row. Similarly, the adjacency matrix corresponding to $I(\alpha_0, \alpha_1)$ is also asymptotically 0-sparse as there are at most $O(\log(n))$ terms on each of its row. Finally, we need to check that the thresholded $\Sigma_p$ is a good approximation of $\Sigma_p$. Recall that for a real symmetric matrix $M$, $\|M\|_2 \leq \max_i(\sum_j |m_{i,j}|)$. (See, e.g., [7] or [24], page 70.) Now, $\sum_{k \geq k_0} \rho^k = \rho^{k_0}/(1-\rho)$, so $\|\Sigma_p - T_{\alpha_1}(\Sigma_p)\|_2 \leq n^{-\alpha_1}/(1-|\rho|)$, which tends to 0 as $n$ goes to infinity. So we conclude that Proposition 1 applies and thresholding the sample covariance (resp., correlation) matrix corresponding to this population covariance will yield an operator norm consistent estimator, a.s., provided the moment conditions are satisfied. In this situation, the moment conditions translate simply into $k \geq 2 + \varepsilon$ for some $\varepsilon$, because $\alpha_0$ can be chosen arbitrarily close to $1/2$ and $\gamma$ arbitrarily close to 0.

Finally, we have the following corollaries that apply to all the theorems and proposition above.

COROLLARY 1 (Infinitely many moments). *Suppose that the entries of $X$ have infinitely many moments. Then all the above results hold with only the sparsity conditions having to be checked.*

COROLLARY 2 (Asymptotic $\beta$-sparsity). *Suppose that the sequence $\Sigma_p$ is asymptotically $\beta$-sparse. Then all the above results apply, with the modification that $\beta$ be replaced by $\widetilde{\beta}_\varepsilon = \beta + \varepsilon$ for $\varepsilon > 0$ but arbitrarily small. In particular, moment conditions need only to be satisfied and checked with $\widetilde{\beta}_\varepsilon$. In the situation where one has infinitely many moments, one therefore only needs to check that the sparsity conditions are satisfied by a $\widetilde{\beta}_\varepsilon$.*

### 3.2.1. *A refinement of Theorem 1.*

We now discuss a refinement of Theorem 1 that allows us to get rid of assumption (ii) there. We remind the reader that this assumption is about the size of the nonzero elements of $\Sigma_p$. The possibility of this refinement was suggested to the author by a question of Professor Peter Bickel whom we thank for his very insightful question. This discussion is included here because it relies on approximation ideas close to the ones we developed for approximating nonsparse matrices by sparse matrices. However, here we will approximate sparse matrices by sparse matrices, the approximating matrix now having quite "large" elements.

Let us first mention the following lemma, which is proved in the Appendix.

LEMMA. *Suppose $M$ is a $p \times p$ real symmetric matrix, which is $\beta$-sparse. Call $m = \max_{i,j} |M_{i,j}|$. Then*

$$\forall k \in 2\mathbb{N} \qquad \|M\|_2 \leq |\operatorname{trace}(M^k)|^{1/k} \leq m p^{\beta(1-1/k)+1/k}(f(k))^{1/k}.$$

We therefore have the following corollary.

COROLLARY 3. *Suppose $\Sigma_p$ is $\beta$-sparse, with $\beta \leq 1/2 - \eta$ and $\eta > 0$. Call $T_{\beta+\varepsilon}(\Sigma_p)$ a version of $\Sigma_p$ thresholded at $Cn^{-(\beta+\varepsilon)}$, where $\varepsilon > 0$ and $C$ is a* *real number (fixed and independent of $p$ and $n$). Call $R_{\beta+\varepsilon} = \Sigma_p - T_{\beta+\varepsilon}(\Sigma_p)$. Assume that $p \asymp n$ as $n \to \infty$. Then*

$$\|\|R_{\beta+\varepsilon}\|\|_2 \to 0 \qquad \text{as } n \to \infty.$$

The conclusion of the previous corollary is that $\beta$-sparse matrices, regardless of the size of their entries, can be approximated in operator norm by $\beta$-sparse matrices whose nonzero elements are greater than $n^{-(\beta+\varepsilon)}$.

PROOF OF COROLLARY 3. We note that $R_{\beta+\varepsilon}(\Sigma_p)$ is $\beta$-sparse, because $\Sigma_p$ is, and the graph corresponding to the adjacency matrix $R_{\beta+\varepsilon}(\Sigma_p)$ is a subgraph of the one corresponding to the adjacency matrix of $\Sigma_p$. Note also that all the entries of $R_{\beta+\varepsilon}(\Sigma_p)$ are less in absolute value than $Cn^{-(\beta+\varepsilon)}$. According to the previous lemma, we therefore have

$$\|\|R_{\beta+\varepsilon}(\Sigma_p)\|\|_2 \leq (f(k))^{1/k} C n^{-(\beta+\varepsilon)} p^{\beta(1-1/k)+1/k}.$$

So if $k \geq 1/\varepsilon$, because our assumption that $p \asymp n$ implies that $p/n$ remains bounded, the right-hand side in the previous equation goes to zero as $n$ goes to infinity. [Recall that $f(k)$ does not depend on $p$.]  □

We are now ready to state our improvement of Theorem 1.

THEOREM 4. *Making the same general assumptions as in Theorems 1, 2 and 3 above [i.e., assumptions* (i)–(iii) *in Theorem 1 are excluded], we now:*

(a) *Assume that $\Sigma_p$ is $\beta$-sparse with $\beta = 1/2 - \eta$ and $\eta > 0$.*

(b) *Pick $\varepsilon_0 > 0$ such that, for some $\delta_0 > 0$, $\beta + \varepsilon_0 \leq 1/2 - \delta_0$.*

(c) *Assume that the random variables $X_{i,j}$ have moments of order $4k$ ($8k$ in the correlation case), with $k$ satisfying the assumptions put forth in Theorem 1, assumption* (iii).

*Then, if $T_{\beta+\varepsilon_0/2}(S_p)$ is the matrix obtained by thresholding the entries of $S_p$ at level $Kn^{-(\beta+\varepsilon_0/2)}$ ($S_p$ having the definition given in Theorems 1, 2 and 3), for some $K > 0$ (fixed and independent of $n$ and $p$), we have*

$$\|\|T_{\beta+\varepsilon_0/2}(S_p) - \Sigma_p\|\|_2 \to 0 \qquad a.s.$$

PROOF. Let us first note that there exists an $\varepsilon_0$ with the characteristics we require. A possible choice is $\varepsilon_0 = \eta/2$, to which $\delta_0 = \eta/2$ could correspond.

The theorem is therefore a consequence of Proposition 1. As a matter of fact, let us pick $\alpha_1 = \beta + \varepsilon_0$ and $\alpha_0 = \beta + \varepsilon_0/4$. As we have seen in Corollary 3, $T_{\beta+\varepsilon_0}(\Sigma_p)$ is $\beta$-sparse and has the property that $\|\|T_{\beta+\varepsilon_0}(\Sigma_p) - \Sigma_p\|\|_2 \to 0$. Clearly, $\alpha_0 < \alpha_1 \leq 1/2 - \delta_0$. In the notation of Proposition 1, $I_{\alpha_1, \alpha_0}$ is a

subset of the set of indices for which $\sigma(i,j) \neq 0$, and hence it is $\beta$-sparse. So in the notation of Proposition 1, our $I_{\alpha_1, \alpha_0}$ is $\gamma$-sparse with $\gamma = \beta \leq \alpha_0 - \zeta_0$, where $\zeta_0 = \varepsilon_0/4$. Note also that the moment assumptions made in Theorem 4 correspond to the moment assumptions made in Proposition 1, with $\gamma = \beta$. So the conclusion of Proposition 1 applies, and in particular, we can take $\alpha = \beta + \varepsilon_0/2$, since $\beta + \varepsilon_0/2 \in (\beta + \varepsilon_0/4, \beta + \varepsilon_0)$. □

3.3. *About 1/2-sparse matrices.* The previous computations are clearly limited to the case where $\beta < 1/2$. A natural question is therefore to ask if this limitation is inherent to the problem, or if it is a consequence of the bounds we use in the mathematical analysis. We now want to highlight the problems that occur in the case $\beta = 1/2$ and show that our result is "sharp": at the level of generality at which we are working, (at least some) 1/2-sparse matrices are not estimable consistently in operator norm by hard thresholding. To show this, we will produce a 1/2-sparse matrix that cannot be consistently estimated in operator norm even at the oracle level. In what follows, we assume that $p/n$ has a finite nonzero limit, $l$, as $n$ tends to infinity.

To do so, we consider estimating a matrix $A$ of the following form:

$$\Sigma_p = \begin{pmatrix} 1 & \alpha_2 & \alpha_3 & \dots & \alpha_p \\ \alpha_2 & 1 & 0 & \dots & 0 \\ \vdots & \vdots & \ddots & \ddots & 0 \\ \alpha_{p-1} & 0 & 0 & 1 & 0 \\ \alpha_p & 0 & 0 & \dots & 1 \end{pmatrix}.$$

To simplify the problem, we assume that the data are multivariate Gaussian, with mean 0, and that we know that the diagonal is composed only of 1's. We estimate $\Sigma_p$ using the sample covariance matrix, putting to 1 the main diagonal, and using the oracle information to put to 0 all other terms except the first row and columns. We call $\widehat{\Sigma}_p$ the corresponding estimator. Note that

$$\Sigma_p - \widehat{\Sigma}_p = \begin{pmatrix} 0 & \alpha_2 - \widehat{\alpha}_2 & \alpha_3 - \widehat{\alpha}_3 & \dots & \alpha_p - \widehat{\alpha}_p \\ \alpha_2 - \widehat{\alpha}_2 & 0 & 0 & \dots & 0 \\ \vdots & \vdots & \ddots & \ddots & 0 \\ \alpha_{p-1} - \widehat{\alpha}_{p-1} & 0 & 0 & 0 & 0 \\ \alpha_p - \widehat{\alpha}_p & 0 & 0 & \dots & 0 \end{pmatrix}.$$

Using the Schur complement formula for determinants (see, e.g., [16], page 22), we see that the characteristic polynomial of this matrix is

$$p(\lambda) = \lambda^{p-2} \left( \lambda^2 - \sum_{i=2}^p (\alpha_i - \widehat{\alpha}_i)^2 \right),$$

and therefore

$$\|\widehat{\Sigma}_p - \Sigma_p\|_2 = \sqrt{\sum_{i=2}^{p} (\alpha_i - \widehat{\alpha}_i)^2}.$$

Note that the computation holds for the corresponding adjacency matrix, giving that $\phi_p(2k) = \mathrm{trace}(A_p^{2k}) = 2(p-1)^k$. So this matrix is 1/2-sparse.

Now, since we assume the data are Gaussian, it is clear that $\lambda_1 = \|\widehat{\Sigma}_p - \Sigma_p\|_2$ has infinitely many moments, using, for instance, Frobenius norm as a bound on $\lambda_1$. Also, $\mathbf{E}(\lambda_1^2) = \sum_{i=2}^{p} \mathbf{E}((\alpha_i - \widehat{\alpha}_i)^2)$. The covariance of elements of the sample covariance matrix is well known in the Wishart case; see, for instance, [2], Theorem 3.4.4, page 87. In our context, we see that $\mathbf{E}((\alpha_i - \widehat{\alpha}_i)^2) = (1 + \alpha_i^2)/(n-1) = \nu_i/(n-1)$. In particular,

$$\mathbf{E}(\lambda_1^2) = \frac{p - 1 + \sum_{i=2}^{p} \alpha_i^2}{n-1} \geq \frac{p-1}{n-1} \to l > 0.$$

We now turn to showing that $\lambda_1^2$ actually converges in probability to this limit.

A standard result in Gaussian multivariate analysis (see [2], Theorem 3.3.2) states that we can write $\widehat{\alpha}_i - \alpha_i = (\sum_{k=1}^{n-1} Z_k)/(n-1)$, where the $Z_k$'s are i.i.d. and mean 0. Hence we get that

$$\mathbf{E}(((\widehat{\alpha}_i - \alpha_i)^2 - \nu_i/(n-1))^2) = \frac{1}{(n-1)^4} \mathbf{E}\left( \sum_{k_1, k_2, k_3, k_4} Z_{k_1} Z_{k_2} Z_{k_3} Z_{k_4} \right).$$

In the above sum, the terms that contain an index repeated only once contribute zero to the expectation. After elementary computations, we see that to first order this expectation is $\mathrm{O}(2\nu_i^2/n^2)$. Using the same ideas (see Appendix), we get that, for $i \neq j$,

$$\mathbf{E}(((\widehat{\alpha}_i - \alpha_i)^2 - \nu_i/(n-1))((\widehat{\alpha}_j - \alpha_j)^2 - \nu_j/(n-1))) = \mathrm{O}\left( \frac{2}{n^2} \alpha_j^2 \alpha_i^2 \vee \frac{1}{n^3} \right).$$

Hence we have that

$$\mathrm{var}(\lambda_1^2) = \mathrm{O}\left( \sum_{i=2}^{p} \frac{2\nu_i^2}{n^2} + \sum_{i \neq j} \frac{2\alpha_j^2 \alpha_i^2}{n^2} \vee \frac{1}{n^3} \right)$$

$$= \mathrm{O}\left( \sum_{i=2}^{p} \frac{2\nu_i^2}{n^2} + \frac{2}{n^2} \left( \sum_{i=2}^{p} \alpha_i^2 \right)^2 \vee \frac{1}{n} \right).$$

Therefore, if, for instance, $\alpha_i = \frac{1}{\sqrt{p}}$, $\mathrm{var}(\lambda_1^2) = \mathrm{O}(\frac{p}{n^2} + \frac{1}{n}) \to 0$ and

$$\lambda_1^2 - \frac{p - 1 + \sum_{i \geq 2} \alpha_i^2}{n} \to 0 \qquad \text{in probability,}$$

and therefore

$$\lambda_1^2 \geq \frac{p-1}{n} \qquad \text{in probability.}$$

Note that if we had tried to estimate $\Sigma_p$ using oracle information about the location of the nonzero coefficient but nothing about the fact the diagonal was equal to 1, we would have encountered the same problem. As a matter of fact, if we call $M_p$ the diagonal matrix with entries $\widehat{\sigma}(i,i)$, we have from previous results in the paper (our moment computations and the 0-sparsity of this matrix) that $\|M_p - \mathrm{Id}_p\|_2 \to 0$ a.s. Note that because $\widehat{\Sigma}_p$ had 1's on its diagonal,

$$\widehat{\Sigma}_p + M_p - \mathrm{Id}_p = \begin{pmatrix} \widehat{\sigma}(1,1) & \widehat{\alpha}_2 & \widehat{\alpha}_3 & \dots & \widehat{\alpha}_p \\ \widehat{\alpha_2} & \widehat{\sigma}(2,2) & 0 & \dots & 0 \\ \vdots & \vdots & \ddots & \ddots & 0 \\ \widehat{\alpha}_{p-1} & 0 & 0 & \widehat{\sigma}(p-1,p-1) & 0 \\ \widehat{\alpha}_p & 0 & 0 & \dots & \widehat{\sigma}(p,p) \end{pmatrix}.$$

So for the oracle estimator that uses only information about the location of the nonzero coefficients, we have

$$\|\widehat{\Sigma}_p + (M_p - \mathrm{Id}_p) - \Sigma_p\|_2 \geq \|\widehat{\Sigma}_p - \Sigma_p\|_2 - \|M_p - \mathrm{Id}_p\|_2 > \frac{p-1}{2n} \qquad \text{a.s.}$$

This example shows that even using oracle information for estimation of the $\Sigma_p$ pointed out above does not lead to an operator norm consistent estimator, in the presence of this simple 1/2-sparse graph. This suggests that for these graphs, simple thresholding might not be a good method. It also suggests that the conditions of our theorems have more to do with the method we propose than with unrefined mathematical details in its analysis.

3.3.1. *Complement on this example.* In what follows we investigate in more details the case where $\alpha_i = 1/\sqrt{p}$. One might ask whether, despite the fact that $\|\Sigma_p - \widehat{\Sigma}_p\|_2$ does not go to zero, $\widehat{\Sigma}_p$ does not have some good characteristics as an estimator of $\Sigma_p$ anyway. In what follows, we show that for both the eigenvalues and eigenvectors, this is not the case.

The previous computations essentially show that

$$\mathbf{E}(\lambda_1^2(\widehat{\Sigma}_p - \mathrm{Id}_p)) = \sum_{i \geq 2} \alpha_i^2 + \frac{1+\alpha_i^2}{n-1} = (\lambda_1(\Sigma_p - \mathrm{Id}))^2 + \frac{p-1}{n-1} + \frac{p-1}{p(n-1)},$$

so at the level of eigenvalues, the answer is negative. Note that the eigenvectors of $\Sigma_p - \mathrm{Id}_p$ and therefore of $\Sigma_p$ are known. The ones corresponding to the nonzero eigenvalues are, calling $\lambda_+ = \sqrt{\sum_{i \geq 2} \alpha_i^2}$,

$$u_+ = \frac{1}{\sqrt{2}\lambda_+} \begin{pmatrix} \lambda_+ \\ \alpha_2 \\ \vdots \\ \alpha_p \end{pmatrix} \quad \text{and} \quad u_- = \frac{1}{\sqrt{2}\lambda_+} \begin{pmatrix} \lambda_+ \\ -\alpha_2 \\ \vdots \\ -\alpha_p \end{pmatrix}.$$

We call $\widehat{u}_+$ the eigenvector corresponding to the positive eigenvalue of $\widehat{\Sigma}_p - \mathrm{Id}_p$. When $\alpha_i = 1/\sqrt{p}$, $\mathrm{cov}(\widehat{\alpha}_i, \widehat{\alpha}_j) = (1_{i=j} + 1/p)/(n-1)$ and $\lambda_+ = \sqrt{(p-1)/p}$. Using the expression above for the eigenvectors, we have

$$2\lambda_+ \widehat{\lambda}_+ \langle u_+, \widehat{u}_+ \rangle = \lambda_+ \widehat{\lambda}_+ + \frac{1}{\sqrt{p}} \sum_{i \geq 2} \widehat{\alpha}_i.$$

Now $\mathrm{var}(\sum_{i \geq 2} \widehat{\alpha}_i) = (p-1)(1+1/p)/(n-1) + (p-1)(p-2)/(p(n-1))$, and $\mathbf{E}(\sum_{i \geq 2} \widehat{\alpha}_i) = (p-1)/\sqrt{p}$, from which we conclude that

$$\frac{1}{\sqrt{p}} \left( \sum_{i \geq 2} \widehat{\alpha}_i \right) - \left( 1 - \frac{1}{p} \right) \to 0 \qquad \text{in probability.}$$

Since when all $\alpha_i = 1/\sqrt{p}$,

$$\mathrm{var}\left( \sum_{i \geq 2} \widehat{\alpha}_i^2 \right) \leq 2 \left( \mathrm{var}\left( \sum_{i \geq 2} (\widehat{\alpha}_i - \alpha_i)^2 \right) + \frac{4}{p} \mathrm{var}\left( \sum_{i \geq 2} \widehat{\alpha}_i \right) \right),$$

the above computations show that, since $p/n \to l$, $\widehat{\lambda}_+ \to \sqrt{1+l}$ in probability and therefore, using Slutsky's lemma, we get that

$$\langle u_+, \widehat{u}_+ \rangle \to \frac{1}{2} \left( 1 + \frac{1}{\sqrt{1+l}} \right) \qquad \text{in probability.}$$

So when $p/n$ has a finite nonzero limit, the angle between these two vectors has a finite nonzero limit (in probability), showing that the eigenvectors are not consistently estimated.

3.4. *Discussion.* In the following, we call $\widehat{\Sigma}_p$ our (final) estimator of $\Sigma_p$, which is obtained from the standard estimator $S_p$. As above, we denote $\Delta_p = \widehat{\Sigma}_p - \Sigma_p$, $\Xi_p = \mathrm{oracle}(S_p) - \Sigma_p$, where $\mathrm{oracle}(S_p)$ is the oracle version of $S_p$, and $D_p = S_p - \Sigma_p$.

3.4.1. *Finite-dimensional character and sharpening of the bounds.* As is clear from the proofs, all the bounds we derive are valid at $n$ and $p$ fixed. Essentially, we get bounds on the probability of deviation of the largest eigenvalue of the matrix $\Delta_p$ from 0. These bounds are polynomial in nature since we used Chebyshev's inequality and worked with moments.

Note that in particular cases, such as when the entries of the data matrix are bounded or satisfy certain tail conditions, these bounds can be sharpened by using (exponential or Gaussian) concentration inequalities for the difference $d_{i,j} = \hat{\sigma}(i,j) - \sigma(i,j)$. If the entries of $X$ are bounded in absolute value by a constant $C$, in the setting of Theorem 1, Hoeffding's inequality (see [15]) would, for instance, give that

$$P(|d_{i,j}| > t) = P(|\hat{\sigma}(i,j) - \sigma(i,j)| > t) \leq 2\exp(-nt^2/(2C^4)).$$

This is a simple consequence of the fact that $\hat{\sigma}(i,j)$ is a sum of i.i.d. random variables and their mean is $\sigma(i,j)$. (Of course, a slight adjustment is needed when dealing with sample covariance matrices, but it does not change the exponential character of the bounds. We give the argument in the simplest case where $S_p = X^*X/n$, the Gaussian MLE when we know the mean is zero.) Suppose that the nonzero coefficients of $\Sigma_p$ are bounded below, in absolute value by $C_1 n^{-1/2+b}$. If we call $B_p$ the event $B_p = \{$at least one mistake is made by the thresholding procedure$\}$, and if we decide to refine our thresholding to a $(\log(n))^a/\sqrt{n}$ threshold, we see, using a simple union bound, that

$$P(B_p) \leq 2p^2(\exp(-(\log n)^{2a}/(2C^4)) + \exp(-((\log n)^a - C_1 n^b)^2/(2C^4))).$$

Therefore, by adding assumptions to our problem, we are able to get sharper bounds on the probability of making a mistake by thresholding.

We can also get better bounds on the probability that $\|\!|\Xi_p|\!\|_2 > \varepsilon$ and $\|\!|\Delta_p|\!\|_2 > \varepsilon$. We assume that $\Sigma_p$ is $\beta$-sparse and use the corresponding notation. Of course, the event $\|\!|\Xi_p|\!\|_2 > \varepsilon$ is contained in the event $\mathrm{trace}(\Xi_p^{2k}) > \varepsilon^{2k}$, which is contained in the event $\max |w_\gamma(2k)| > \varepsilon^{2k}/(f(k)p^{1+\beta(2k-1)})$ which is contained in the event $\max |d_{i,j}| > \varepsilon/(f(k)^{1/(2k)}p^{1/2k+\beta(1-1/2k)})$. Hence, by using Hoeffding's inequality, we get

$$P(\|\!|\Xi_p|\!\|_2 > \varepsilon) \leq 2p^2 \exp(-n\varepsilon^2 p^{-2\beta} p^{(\beta-1)/k}/(2C^4 f(k)^{1/k})).$$

Finally, using the fact that $\{\|\!|\Delta_p|\!\|_2 > \varepsilon\} \subseteq (\{\|\!|\Xi_p|\!\|_2 > \varepsilon\} \cap B_p^c) \cup B_p$, we see that

$$P(\|\!|\Delta_p|\!\|_2 > \varepsilon) \leq P(B_p) + P(\|\!|\Xi_p|\!\|_2 > \varepsilon),$$

for which we just derived bounds. Similar types of bounds can be obtained in the context of Theorem 2, when, for instance, Hoeffding's inequality applies.

Though these results are sharper than the ones announced in the theorems above, they are less general. Because one of our concerns was distributional generality, we decided to give the theorems in general form with less sharp bounds.

3.4.2. *Beyond the bounded $p/n$ assumption.* A close look at the proofs of the theorems and the bounds above reveals that the assumption that $p/n$ remains bounded can be relaxed. As a matter of fact, our bounds on expected values of traces are generically of the form $\mathrm{O}(p^\gamma n^{-\lambda})$, and all we require is that this quantity goes to zero fast enough. If we focus on the oracle version of the theorems, we see that the bounds are of the form

$$\mathbf{E}(\mathrm{trace}(\Xi_p^{2k})) = \mathrm{O}(n^{-k}p^{1+\beta(2k-1)}).$$

If $p = \mathrm{O}(n^\nu)$, we see that the exponent in $n$ becomes of the form $k(2\beta\nu - 1) + \nu(1-\beta)$. If this quantity is less than $-(1+\varepsilon)$ for some $\varepsilon > 0$ and $k = k_0$, then we will have a.s. convergence of $\Xi_p$ to zero in operator norm. This condition is satisfied if

$$\nu \le \frac{k - (1+\varepsilon)}{1 + \beta(2k-1)}.$$

So in particular, if we are working with random variables with infinitely many moments, the oracle results will hold almost surely for a $\beta$-sparse matrix when

$$p = \mathrm{O}(n^{1/(2\beta)-\eta}) \qquad \text{for some } \eta \text{ arbitrarily small.}$$

As a matter of fact, all we need to do is pick a finite number $k_1$ such that

$$\frac{k_1 - (1+\varepsilon)}{1 + \beta(2k_1 - 1)} > 1/(2\beta) - \eta$$

and carry out the analysis for $\mathbf{E}(\mathrm{trace}(\Xi_p^{2k_1}))$. $k_1$ exists (and is finite), because $\frac{k-(1+\varepsilon)}{1+\beta(2k-1)} \to 1/(2\beta)$, as $k$ goes to infinity. If there are only $4k_1$ moments, the results will hold, too.

On the other hand, the nonoracle results will be satisfied in the context of Theorem 1 as soon as

$$\nu \le \frac{k(1-2\alpha_0) - (1+\varepsilon)}{2},$$

a constraint much less restrictive than the previous one in general. Finally, we note that Proposition 1 would apply if, assuming the other constraints were satisfied, we also had

$$\nu \le \frac{k - (1+\varepsilon)}{1 + \gamma(2k-1)}.$$

### 3.5. *Consequences of spectral norm convergence.*

#### 3.5.1. *Convergence of eigenvalues.*

We recall some classical facts from matrix analysis. Firstly, if $A$ and $B$ are two symmetric matrices, and if $\lambda_i$ is their $i$th eigenvalue, where the eigenvalues are sorted in decreasing order, we have, by Weyl's theorem (Theorem 4.3.1 in [16])

$$|\lambda_i(A) - \lambda_i(B)| \leq \|A - B\|_2.$$

Because the matrix $S_p$ is symmetric, the thresholded version of it is symmetric, too. Therefore the operator norm convergence we showed implies the following:

FACT 4. *When the thresholded estimator $\widehat{\Sigma}_p$ is a spectral norm consistent estimator of the population covariance or correlation matrix $\Sigma_p$, all the eigenvalues of $\widehat{\Sigma}_p$ are consistent estimators of the population eigenvalues.*

#### 3.5.2. *Convergence of eigenvectors.*

Perhaps even more interestingly, controlling the spectral norm allows us to get very good control on the angles between the eigenspaces of the population and sample covariance matrix, through the use of the classical $\sin(\theta)$ theorems of Davis and Kahan ([10], Section 2, and [24], Section V.3). For the sake of completeness we quote a version of this important result (Theorem 2 in [10]) and show how to exploit it in our context.

THEOREM 5 [$\sin(\theta)$ theorem]. *Suppose $\Sigma_p$ has the spectral resolution*

$$\begin{pmatrix} X_1' \\ X_2' \end{pmatrix} \Sigma_p (X_1 X_2) = \operatorname{diag}(L_1, L_2)$$

*with $(X_1 X_2)$ an orthogonal matrix, $X_1$ being a $p \times k$ matrix. Suppose $Z$ is a $p \times k$ matrix with orthogonal columns, and for any Hermitian matrix $M$ of order $k$, call $R = \Sigma_p Z - ZM$. Suppose the eigenvalues of $M$ are contained in an interval $[\alpha, \beta]$ and that for some $\delta > 0$, the eigenvalues of $L_2$ are contained in $\mathbb{R} \setminus [\alpha - \delta, \beta + \delta]$. Then for any unitarily invariant norm,*

$$\|\sin \Theta[\mathcal{R}(X_1), \mathcal{R}(Z)]\| \leq \frac{\|R\|}{\delta},$$

*where $\Theta[\mathcal{R}(X_1), \mathcal{R}(Z)]$ stands for the canonical angles between the column space of $X_1$ and that of $Z$, and $\sin \Theta[\mathcal{R}(X_1), \mathcal{R}(Z)]$ is the diagonal matrix containing these angles.*

These angles are closely connected to canonical correlation analysis: their cosines are the canonical correlations for the "data matrices" $X_1$ and $Z$.

We therefore have the following corollary to Theorems 2 and 3:

COROLLARY 4 (Consistency of eigenspaces). *Suppose $\Sigma_p$ has a group of eigenvalues contained in an interval and separated from the other eigenvalues by $\delta > 0$. Call the set of their indices (after, say, ordering them) $J$. Then the canonical angles between the column space of the corresponding eigenvectors and the column space of the eigenvectors of $\widehat{\Sigma}_p$ (our thresholding estimator) corresponding to the eigenvalues of $\widehat{\Sigma}_p$ with index set $J$ goes to zero a.s.*

PROOF. Call $\widehat{\lambda}_j$ the eigenvalues of $\widehat{\Sigma}_p$ with index set $J$. Let $M$ be the diagonal matrix with diagonal entries the $\{\widehat{\lambda}_j\}$. Call $L_2$ the set consisting of the other eigenvalues of $\Sigma_p$. Note that the convergence of eigenvalues guarantees that the $\{\widehat{\lambda}_j\}_{j \in J}$ will a.s. stay away from $L_2$, by a distance at least $\delta_2 > 0$. Call $Z_j$ the eigenvectors corresponding to $\widehat{\lambda}_j$ and $Z$ the matrix with columns $Z_j$ (if some eigenvalues have multiplicity higher than 1, then we pick a set of such eigenvectors). We can write $\Sigma_p = \widehat{\Sigma}_p - \Delta_p$ with $\|\!|\Delta_p\|\!|_2 \to 0$ a.s. Note that $\widehat{\Sigma}_p Z = ZM$, so $\Sigma_p Z = ZM - \Delta_p Z$. Therefore $R = -\Delta_p Z$ and because $\|\!| \cdot \|\!|_2$ is matrix norm and the columns of $Z$ are orthonormal, $\|\!|R\|\!|_2 \leq \|\!|\Delta_p\|\!|_2$. Applying Theorem 5 with these inputs gives the result. □

3.6. *Practical considerations.* The theoretical part of this paper essentially says $\beta$-sparse matrices with $\beta < 1/2$ are asymptotically estimable, in the strong notion of estimability induced by the spectral norm. However, it does not give much information about how to choose the thresholding parameter.

In practice, covariance matrices are estimated for a purpose other than simply estimating them. So in concrete applications, users would most likely be able to find a penalty function that incorporates a measure of performance of a certain estimator and mitigates it with how sparse the corresponding matrix is. Then cross-validation or resampling techniques might be used to assess the performance of different estimators and choose the threshold from the data. Note also, that in [7], Section 5, the authors propose a technique for choosing a banding parameter from the data, which is shown empirically to work quite well. Such technique is transferable in our context, through some fairly straightforward steps.

However, a shortcoming of resampling techniques is their heavy computational cost. Thresholding methods are appealing because they are easily "parallelizable" and can be used on very large dimensional datasets. Therefore having an a priori method that works reasonably well and is not too computationally expensive is also worthwhile. Of course there is a clear link between thresholding and testing the hypothesis that a certain parameter is 0. As a practical ansatz, one method that can be tried is the following: get a $p$-value for the hypothesis $\sigma(i, j) = 0$ for all $i > j$. Such a $p$-value can be obtained by bootstrap methods and since we are dealing with means those reduce to a simple $z$-test. Then perform the Benjamini–Hochberg procedure (see [5]) for these $p$-values, using the FDR parameter $1/\sqrt{p}$. Though the theoretical part of [5] does not apply, we found in the practical examples we ran (limited to Gaussian simulations and relatively simple population covariance matrices) that this worked reasonably well. We include some figures illustrating our simulations (see Appendix A.2). If speed is the most important issue, not using the FDR but testing each entry at level $\alpha/\sqrt{p}$ seems also to yield reasonable results.

*The issue of positive semidefiniteness.* We note that it is possible that our estimators will not be positive definite: thresholding entrywise the sample covariance or correlation matrix does not guarantee positive definiteness of the resulting estimator. Our theorems, however, say that if the population matrices have a smallest eigenvalue bounded away from zero (uniformly in $p$), then asymptotically our estimators will yield positive definite matrices (in that case, the theorems also imply spectral norm consistency of $\widehat{\Sigma}_p^{-1}$ for $\Sigma_p^{-1}$). If, in practice, one encounters a nonpositive definite estimator, it is clear that the problem at hand should dictate the strategy to remedy this flaw. Three general ideas can nevertheless be applied: one might think of "projecting" the estimator on the cone of positive semidefinite matrices, using semi definite programming and probably a sparseness penalty. The feasibility of this idea depends of course on the dimensionality of the problem and it is unlikely to work well (at this point in time) in truly high dimension. Another idea would be to do a singular value decomposition of the estimator, which is possible even in high dimension, since the estimator is by construction sparse, and hence falls within the reach of several fast algorithms in numerical linear algebra. Then one could keep a smaller rank approximation of $\widehat{\Sigma}_p$ as the final estimator, $\widehat{\Sigma}_f$, by putting, for instance, all the negative eigenvalues of $\widehat{\Sigma}_p$ to zero [or instead of 0 a real $g(p)$, with $g(p) \to 0$]. Note that $\widehat{\Sigma}_f$ can also be shown to be a consistent estimator of the population covariance, in spectral norm, since $\|\widehat{\Sigma}_f - \widehat{\Sigma}_p\|_2 \to 0$ because the negative eigenvalues of $\widehat{\Sigma}_p$ have to converge to zero (otherwise $\|\widehat{\Sigma}_p - \Sigma_p\|_2$ would not tend to 0). The main drawback of such a solution to the positive definiteness problem is that we may lose the sparsity of the estimator, a feature that is in general desirable. However, its spectral characteristics would be quite easy to obtain, even in high dimension. The third idea would be to consider for our estimator, in the case where $\widehat{\Sigma}_p$ turns out to not be positive definite, the matrix $\widehat{\Sigma}_f = \widehat{\Sigma}_p - \lambda_p(\widehat{\Sigma}_p)\mathrm{Id}_p$, where $\lambda_p(\widehat{\Sigma}_p)$ is the smallest eigenvalue of $\widehat{\Sigma}_p$. Since, as we just noted, $|\lambda_p(\widehat{\Sigma}_p)| \to 0$, we see that $\|\widehat{\Sigma}_f - \widehat{\Sigma}_p\|_2 \to 0$ and hence $\|\widehat{\Sigma}_f - \Sigma_p\|_2$ tends to 0, too. Hence, $\widehat{\Sigma}_f$ is operator norm consistent. Note that it is also sparse, because adding a diagonal

matrix does not change anything about the sparsity of our estimator. So this is a sparse positive semidefinite and operator norm consistent estimator of $\Sigma_p$.

*Robustness issues.* Finally, we note that the results of this paper suggest that acting entrywise on the sample covariance matrix is a way to create good estimators of $\Sigma_p$. In particular, when other issues such as robustness or contamination by heavy-tailed data arise, using (entrywise) more robust estimators than the sample covariance is likely to give improved results.

**4. Conclusion.** In this paper we have investigated the theoretical properties of the idea of thresholding the entries of a sample covariance (or correlation) matrix to better estimate the population covariance, when it is assumed (or known) to be sparse. We have shown that the natural notion of sparsity, coming from problems concerning random vectors, is not appropriate when one is concerned with estimating matrices and in particular their spectral properties. By contrast, we propose an alternative notion of sparsity, based on properties of the graph corresponding to the adjacency matrix of the population covariance. We have shown that our notion of sparsity divides sharply classes of matrices that are estimable through hard thresholding and those that are not, an appealing property. The notion of sparsity we propose is invariant under permutation of the order of the variables and hence is well suited for the analysis of problems where there is no canonical ordering of the variables. It is also related to the spectral norm of the adjacency matrix of the population covariance.

We show that $\beta$-sparse matrices, with $\beta < 1/2$, are consistently estimable in operator (a.k.a. spectral) norm, a strong notion of convergence that implies consistency of all eigenvalues and eigenspaces corresponding to eigenvalues separated from the rest of the spectrum (see Section 3.5). Practically, the results of simulations are maybe not as striking as one may have hoped for, but lead to great improvement over the sample covariance (or correlation) matrix.

We also show that certain nonsparse matrices are estimable by sparse matrices through the thresholding method we analyzed. Numerically, this method has many advantages in terms of implementation. It is easy to implement, and leads to sparse matrices, which have the desirable property that their eigenvalues and eigenvectors can be numerically computed efficiently, even in high dimension. Also, since the method acts in an entrywise fashion, the corresponding algorithm is easily parallelizable and in general produces results quickly.

Statistically, our results mean that under the assumption of $\beta$-sparsity, $\beta < 1/2$, applying the natural practical idea of thresholding the entries of a sparse matrix leads to good theoretical convergence properties. However, we also show that in situations that are not inconceivable in practice, that is, $\beta \geq 1/2$, this strategy may sometimes fail to give an estimator as good as what we required. More sophisticated approaches may be needed in these more difficult cases, though, as noted above, the simple thresholding approach has even then many practical virtues.

## APPENDIX

**A.1. 1/2-sparse matrices: details of computations.** In what follows, we use the notation $N$ for the quantity $n - 1$ (so $N = n - 1$) in an effort to alleviate the notation. The computations that follow are used in Section 3.3 and the notation are defined there. Recall that $\nu_i = 1 + \alpha_i^2$ and $i \geq 2$. We give a detailed explanation of our estimate of

$$\mathbf{E}(((\widehat{\alpha}_i - \alpha_i)^2 - \nu_i/N)((\widehat{\alpha}_j - \alpha_j)^2 - \nu_j/N)).$$

Clearly, the only thing we need to control is $\mathbf{E}((\widehat{\alpha}_i - \alpha_i)^2(\widehat{\alpha}_j - \alpha_j)^2)$, since $\nu_i/N$ and $\nu_j/N$ are the means of $(\widehat{\alpha}_i - \alpha_i)^2$ and $(\widehat{\alpha}_j - \alpha_j)^2$. Note that we can write $(\widehat{\alpha}_i - \alpha_i) = \sum_{k=1}^{N} Z_k(i)/N$, where the $Z_k(i)$'s are i.i.d. and mean 0. Similarly, we can write $(\widehat{\alpha}_j - \alpha_j) = \sum_{k=1}^{N} Y_k(j)/N$. Note that $Y_k(j)$ is independent of $Z_l$ if $k$ is different from $l$. Therefore,

$$\mathbf{E}((\widehat{\alpha}_i - \alpha_i)^2(\widehat{\alpha}_j - \alpha_j)^2) = \frac{1}{N^4}\mathbf{E}\Big(\sum Z_{k_1}(i)Z_{k_2}(i)Y_{k_3}(j)Y_{k_4}(j)\Big).$$

In the previous sum if an index appears only once in the product, the expectation is zero. So only terms where each index appears an even number of times will matter.

We first focus on terms where we have two distinct indices; the contribution of such terms is

$$\frac{N(N-1)}{N^4}\mathbf{E}(Z_1^2 Y_2^2 + Z_1 Z_2 Y_1 Y_2 + Z_1 Z_2 Y_2 Y_1).$$

We can limit our investigations to the terms with two distinct indices since there are only $N$ terms of the form $Z_1^2 Y_1^2$, so their contribution will be asymptotically negligible. Now, $\mathbf{E}(Z_1^2 Y_2^2) = \nu_i \nu_j$, by independence and definition. Also, if $X$ is multivariate Gaussian vector with covariance $\Sigma_p$,

$$\mathbf{E}(Z_1(i)Y_1(j)) = \mathbf{E}((X_1 X_i - \alpha_i)(X_1 X_j - \alpha_j)) = \mathbf{E}(X_1^2 X_i X_j - \alpha_i \alpha_j)$$

$$= \sigma(1,1)\sigma(i,j) + \alpha_i \alpha_j + \alpha_j \alpha_i - \alpha_i \alpha_j$$

$$= \alpha_i \alpha_j + 1_{i=j},$$

by using the fact that we are working with Gaussian random variables. Therefore, if $i \neq j$,

$$\mathbf{E}\left(\left\{(\widehat{\alpha}_i - \alpha_i)^2 - \frac{\nu_i}{N}\right\}\left\{(\widehat{\alpha}_j - \alpha_j)^2 - \frac{\nu_j}{N}\right\}\right)$$

$$= \left(\frac{1}{N^2} - \frac{1}{N^3}\right)[2(\alpha_i\alpha_j)^2 + \nu_i\nu_j] + \frac{\mathbf{E}(Z_1^2 Y_1^2)}{N^3} - \frac{\nu_i\nu_j}{N^2},$$

$$= 2(\alpha_i\alpha_j)^2 \left(\frac{1}{N^2} - \frac{1}{N^3}\right) + \frac{\mathbf{E}(Z_1^2 Y_1^2) - \nu_i\nu_j}{N^3},$$

$$= \mathrm{O}\left(\frac{(\alpha_i\alpha_j)^2}{N^2} \vee \frac{1}{N^3}\right).$$

In the case where $\alpha_i\alpha_j = \mathrm{o}(1/\sqrt{N})$, we see that this term is of order $1/n^3$.

**A.2. Performance of estimator: graphical illustration.** The images of this subsection illustrate the performance of the estimator, assessing visually its variability and comparing it to the sample covariance matrix. All simulations were done with Gaussian data; the thresholding was made according to the FDR rule—in connection with $z$-tests—with FDR parameter $1/\sqrt{p}$. Our illustrations focus on the properties of eigenvalues because they are easier to visualize.

All matrices investigated are (symmetric) Toeplitz matrices, because of the ease with which they can be simulated. We did not randomly permute the "variables" because this would have had no effect on the performance of the estimator; in particular, the eigenvalues would be exactly the same. These matrices can be summarized by their first row, which is what we refer to when speaking of "coefficients" below.

*Case of a Toeplitz matrix, with $n = p = 500$, and coefficients $(1, 0.3, 0.4,$ $0, \ldots, 0)$.* This situation should be fairly easy since the nonzero coefficients are quite large compared to the variance of $\widehat{\sigma}(i,j)$'s for those $(i,j)$ for which $\sigma(i,j) = 0$. The results are illustrated in Figure 1(a).

*Case of a Toeplitz matrix, with $n = p = 500$, and coefficients $(2, 0.2, 0.3,$ $0, -0.4, 0, \ldots, 0)$.* This situation is a bit harder than the one above a priori, as the nonzero coefficients are not as large compared to the variance of $\widehat{\sigma}(i,j)$'s for those $(i,j)$ for which $\sigma(i,j) = 0$ as they are in the previous example. The results are illustrated in Figure 2(a).

*Case of a nonsparse Toeplitz matrix, with $n = 500$, $p = 100$, and coefficients $\{0.3^k\}_{k=0}^{p-1}$.* This situation illustrates the approximation of a nonsparse matrix by a sparse matrix. As seen above, this population covariance can be approximated in spectral norm by a 0-sparse matrix. In these types of situations, it is possible that thresholding might be a bit "harsh" and "smoother" regularization approaches might lead to better empirical results. The results are illustrated in Figure 3(a).

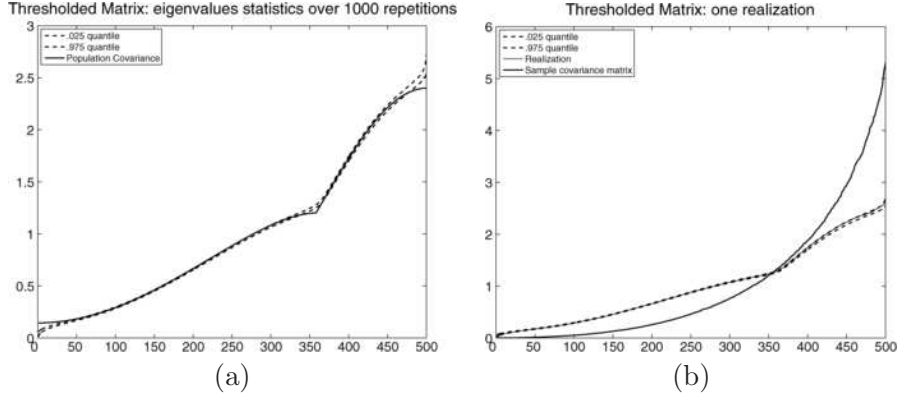

FIG. 1. *Case of a Toeplitz* $(1, 0.3, 0.4, 0, \ldots, 0)$ *population covariance matrix* $\Sigma_p$, $n = p = 500$. *The dashed lines correspond to the 0.025 and 0.975 quantiles of the empirical distribution of the kth eigenvalue, for $k = 1$ to $p$. The data were $\mathcal{N}(0, \Sigma_p)$ and the experiment was repeated 1000 times. As we can see, the estimator is very stable. It does well, especially "far" from the edges of the spectrum. For this particular $\Sigma_p$, it can be explained by the fact that the nonzero coefficients in the matrix are easily detectable, when $n = 500$. The improvement over the sample covariance matrix is quite dramatic.* (a) *Variability of estimator and population spectrum: scree plot of population and corresponding confidence bounds for ordered eigenvalues of our estimator.* (b) *Comparison between scree plot of our estimator (a.k.a. "Realization": the continuous line between the two dashed ones) and that of the sample covariance matrix on one realization, picked at random from our 1000 repetitions.*

**A.3. Some linear algebraic results.** In the course of our proofs, we need the following two lemmas, which are also of independent interest.

We first prove the following lemma, which we needed earlier in the paper.

LEMMA A.1. *Suppose $M$ is a $p \times p$ real symmetric matrix, which is $\beta$-sparse. Call $m = \max_{i,j} |M_{i,j}|$. Then*

$$\forall k \in 2\mathbb{N} \qquad \|M\|_2 \leq |\operatorname{trace}(M^k)|^{1/k} \leq m p^{\beta(1 - 1/k) + 1/k} (f(k))^{1/k}.$$

PROOF. In what follows, $k$ is an even integer. Let $\gamma$ be a closed walk of length $k$. Then, $w_\gamma$, its weight, clearly satisfies, according to Definition 2,

$$|w_\gamma| \leq m^k.$$

So, since

$$\operatorname{trace}(M^k) = \sum_{\gamma \in \mathcal{C}_p(k)} w_\gamma,$$

we clearly have

$$|\operatorname{trace}(M^k)| \leq m^k \phi_p(k).$$

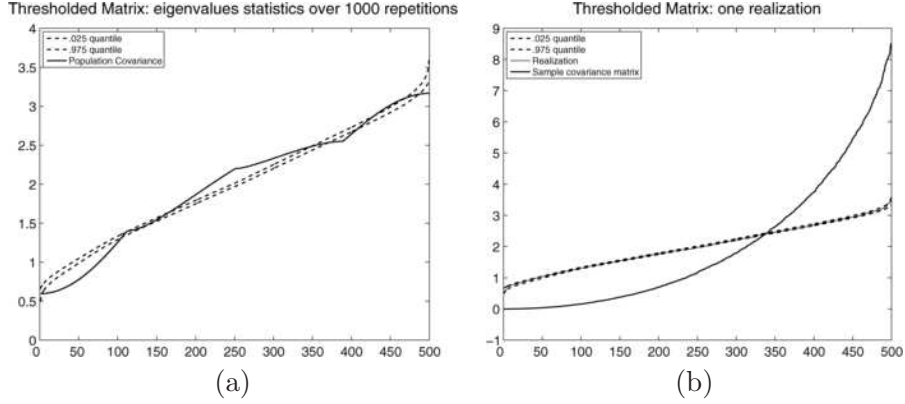

FIG. 2. *Case of a Toeplitz* $(2, 0.2, 0.3, 0, -0.4, 0, \ldots, 0)$ *population covariance matrix* $\Sigma_p$, $n = p = 500$. *The dashed lines correspond to the 0.025 and 0.975 quantiles of the empirical distribution of the kth eigenvalue, for* $k = 1$ *to* $p$. *The data were* $\mathcal{N}(0, \Sigma_p)$ *and the experiment was repeated 1000 times. As we can see, the estimator is very stable. It does capture the support of the spectrum fairly accurately, but is not as good in capturing the fine details of the bulk. For this particular* $\Sigma_p$, *there is (compared to the previous example of Figure 1) a certain lack of accuracy when estimating the adjacency matrix* $A_p$ *of* $\Sigma_p$, *when* $n = 500$. *The improvement over the sample covariance matrix is quite dramatic.* (a) *Variability of estimator and population spectrum: scree plot of population and corresponding confidence bounds for ordered eigenvalues of our estimator.* (b) *Comparison between scree plot of our estimator (a.k.a. "Realization": the continuous line between the two dashed ones) and that of the sample covariance matrix on one realization, picked at random from our 1000 repetitions.*

Since we assume that $M$ is $\beta$-sparse,

$$|\operatorname{trace}(M^k)| \leq f(k) p^{\beta(k-1)+1} m^k. \qquad \square$$

We now turn to another result we needed in the course of our proofs.

LEMMA A.2. *Suppose that* $A$ *and* $B$ *are two real symmetric* $p \times p$ *matrices, with* $|A(i, j)| \leq B(i, j)$. *Then,*

$$\|A\|_2 \leq \|B\|_2.$$

PROOF. Recall that, in the notation of Definition 2,

$$\operatorname{trace}(A^k) = \sum_{\gamma \in \mathcal{C}_p(k)} w_\gamma(A).$$

Now, we clearly have, if $\gamma$ is the walk $i_1 \to i_2 \to \cdots \to i_k \to i_{k+1} = i_1$,

$$|w_\gamma(A)| = |A(i_1, i_2) \cdots A(i_k, i_{k+1})|$$
$$\leq |A(i_1, i_2)| \cdots |A(i_k, i_{k+1})|$$
$$\leq B(i_1, i_2) \cdots B(i_k, i_{k+1}) = w_\gamma(B).$$

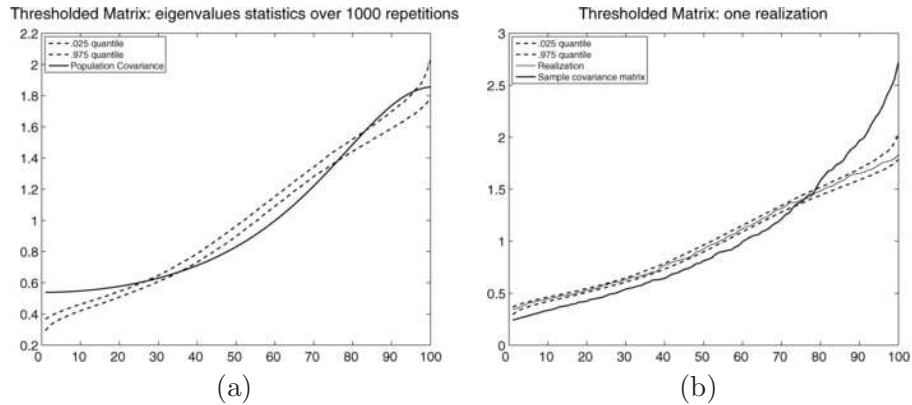

FIG. 3. *Case of a Toeplitz $\{0.3^k\}_{k=0}^{p-1}$ population covariance matrix $\Sigma_p$, $n = 500, p = 100$. The dashed lines correspond to the 0.025 and 0.975 quantiles of the empirical distribution of the kth eigenvalue, for $k = 1$ to p. The data were $\mathcal{N}(0, \Sigma_p)$ and the experiment was repeated 1000 times. As we can see, the estimator is very stable. The problem is harder for the thresholding technique than the one illustrated in Figure 1, and it is possible that less "harsh" regularizations might perform slightly better. The improvement over the sample covariance matrix is still quite dramatic.* (a) *Variability of estimator and population spectrum: scree plot of population and corresponding confidence bounds for ordered eigenvalues of our estimator.* (b) *Comparison between scree plot of our estimator (a.k.a. "Realization": the continuous line between the two dashed ones) and that of the sample covariance matrix on one realization, picked at random.*

So, if $|A|_{\text{Had}}$ is the matrix with $(i, j)$ entry $|A(i, j)|$, we have

$$|\operatorname{trace}(A^k)| \leq \operatorname{trace}(|A|_{\text{Had}}^k) \leq \operatorname{trace}(B^k).$$

For real symmetric matrices, we have $\|A\|_2 = \lim_{k \to \infty} [\operatorname{trace}(A^{2k})]^{1/(2k)}$, and therefore we can conclude that

$$\|A\|_2 \leq \||A|_{\text{Had}}\|_2 \leq \|B\|_2. \qquad \square$$

**Acknowledgments.** The author is very grateful to Peter Bickel for many very interesting discussions on this and related topics. In particular, Theorem 4, a refinement of the main theorems of the paper, came out of a very insightful question posed by Professor Bickel to the author. The author would also like to thank Elizabeth Purdom for discussions that led to clarifications at the beginning of this project and Jim Pitman for references.

## REFERENCES

[1] ANDERSON, G. W. and ZEITOUNI, O. (2006). A CLT for a band matrix model. *Probab. Theory Related Fields* **134** 283–338. MR2222385

[2] ANDERSON, T. W. (2003). *An Introduction to Multivariate Statistical Analysis*, 3rd ed. Wiley, Hoboken, NJ. MR1990662

[3] BAI, Z. D. and SILVERSTEIN, J. W. (2004). CLT for linear spectral statistics of large-dimensional sample covariance matrices. *Ann. Probab.* **32** 553–605. MR2040792

[4] BENGTSSON, T. and FURRER, R. (2007). Estimation of high-dimensional prior and posterior covariance matrices in Kalman filter variants. *J. Multivariate Anal.* **98** 227–255. MR2301751

[5] BENJAMINI, Y. and HOCHBERG, Y. (1995). Controlling the false discovery rate: a practical and powerful approach to multiple testing. *J. Roy. Statist. Soc. Ser. B* **57** 289–300. MR1325392

[6] BHATIA, R. (1997). *Matrix Analysis.* Springer, New York. MR1477662

[7] BICKEL, P. J. and LEVINA, E. (2007). Regularized estimation of large covariance matrices. *Ann. Statist.* **36** 199–227. MR2387969

[8] BICKEL, P. J. and LEVINA, E. (2008). Covariance regularization by thresholding. *Ann. Statist.* **36** 2577–2604.

[9] D'ASPREMONT, A., BANERJEE, O. and EL GHAOUI, L. (2008). First-order methods for sparse covariance selection. *SIAM J. Matrix Anal. Appl.* **30** 56–66. MR2399568

[10] DAVIS, C. and KAHAN, W. M. (1970). The rotation of eigenvectors by a perturbation. III. *SIAM J. Numer. Anal.* **7** 1–46. MR0264450

[11] EL KAROUI, N. (2007). Tracy–Widom limit for the largest eigenvalue of a large class of complex sample covariance matrices. *Ann. Probab.* **35** 663–714. MR2308592

[12] EL KAROUI, N. (2008). Spectrum estimation for large dimensional covariance matrices using random matrix theory. *Ann. Statist.* **36** 2757–2790.

[13] GEMAN, S. (1980). A limit theorem for the norm of random matrices. *Ann. Probab.* **8** 252–261. MR0566592

[14] HAFF, L. R. (1980). Empirical Bayes estimation of the multivariate normal covariance matrix. *Ann. Statist.* **8** 586–597. MR0568722

[15] HOEFFDING, W. (1963). Probability inequalities for sums of bounded random variables. *J. Amer. Statist. Assoc.* **58** 13–30. MR0144363

[16] HORN, R. A. and JOHNSON, C. R. (1990). *Matrix Analysis.* Cambridge Univ. Press. MR1084815

[17] HUANG, J. Z., LIU, N., POURAHMADI, M. and LIU, L. (2006). Covariance matrix selection and estimation via penalised normal likelihood. *Biometrika* **93** 85–98. MR2277742

[18] JAMES, W. and STEIN, C. (1961). Estimation with quadratic loss. *Proc. 4th Berkeley Symp. Math. Statist. Probab.* **I** 361–379. Univ. California Press, Berkeley. MR0133191

[19] JONSSON, D. (1982). Some limit theorems for the eigenvalues of a sample covariance matrix. *J. Multivariate Anal.* **12** 1–38. MR0650926

[20] LEDOIT, O. and WOLF, M. (2004). A well-conditioned estimator for large-dimensional covariance matrices. *J. Multivariate Anal.* **88** 365–411. MR2026339

[21] MARČENKO, V. A. and PASTUR, L. A. (1967). Distribution of eigenvalues in certain sets of random matrices. *Mat. Sb. (N.S.)* **72** 507–536. MR0208649

[22] SILVERSTEIN, J. W. (1995). Strong convergence of the empirical distribution of eigenvalues of large-dimensional random matrices. *J. Multivariate Anal.* **55** 331–339. MR1370408

[23] STANLEY, R. P. (1997). *Enumerative Combinatorics.* **I**. Cambridge Univ. Press. MR1442260

[24] STEWART, G. W. and SUN, J. G. (1990). *Matrix Perturbation Theory.* Academic Press, Boston, MA. MR1061154

[25] Wigner, E. (1955). Characteristic vectors of bordered matrices with infinite dimen-
    sions. *Ann. of Math.* **62** 548–564. MR0077805

Department of Statistics
University of California, Berkeley
367 Evans Hall
Berkeley, California 4720-3860
USA
E-mail: nkaroui@stat.berkeley.edu